%% file: ASEE2025.tex
\documentclass{article}

\usepackage[letterpaper,top=2cm,bottom=2cm,left=3cm,right=3cm,marginparwidth=1.75cm]{geometry}

\PassOptionsToPackage{usenames,dvipsnames,svgnames,table,xcdraw}{xcolor}
\PassOptionsToPackage{dvipsnames,svgnames,table,xcdraw}{enumitems}

\input{Defines/definesPackages}

\input{Defines/definesCommands}

\title{Calculus for the Modern Engineer: Putting the Joy Back in Learning Advanced Mathematics}
\author{Jessy Grizzle \\ Department of Robotics, University of Michigan, Ann Arbor, MI 48109}

\begin{document}
\maketitle

\begin{abstract} 
Many engineering students enter college with a passion for math and physics, only to have their excitement dimmed by a rigid, outdated calculus curriculum. The University of Michigan’s Robotics Department is piloting a new course, \textsc{Calculus for the Modern Engineer}, designed to reintroduce the joy of learning advanced mathematics. This 4-credit course integrates Differential and Integral Calculus of a single variable, vector derivatives, and Ordinary Differential Equations (ODEs) into a unified, one-semester curriculum tailored specifically for engineering students. Departing from the traditional calculus sequence codified in the 1950s, the course begins with definite integration---a concept students easily grasp through sums---before advancing through single-sided limits, differentiation, antiderivatives, and ODEs. By leveraging modern computational tools such as Julia, Large Language Models (LLMs), and Wolfram Alpha Pro, the course shifts the focus from tedious hand calculations to conceptual mastery and real-world application. Engineering case studies, covering topics like numerical integration, optimization, and feedback control, allow students to see firsthand how mathematics powers real-world engineering challenges. Supported by an open-source textbook and interactive programming assignments, \textsc{Calculus for the Modern Engineer} aims to rekindle student enthusiasm for mathematics and equip them with essential computational skills for their careers. \href{https://grizzle.robotics.umich.edu/education/rob201.html}{Updates posted here.}
\end{abstract}

\section{Why Calculus Needs a Complete Overhaul in 2024}
\label{sec:Intro}
\input{Sections/01Intro.tex}

\section{Prior and Ongoing Reform Efforts}
\label{sec:reformEfforts}
\input{Sections/02ReformEfforts.tex}

\section{Content of the Open-source Textbook}
\label{sec:Textbook}
\input{Sections/03TextbookContent.tex}

\section{Three Pillar Projects}
\label{sec:Projects}
\input{Sections/04Projects.tex}

\section{Student Assessment and Sentiment}
\label{sec:PreliminaryAssessment}
\input{Sections/05StudentAssessment.tex}

\section{Barriers to Calculus Reform}
\label{sec:Barriers}
\input{Sections/06BarriersToReform.tex}

\section{Discussion on Course Structure and Computational Focus}
\label{sec:Discussion}
\input{Sections/07Discussion.tex}

\section{Conclusions and Future Perspectives}
\label{sec:ConclusionsPerspectives}

\input{Sections/08ConclusionPerspectives}

\section*{Acknowledgment}
\textsc{Calculus for the Modern Engineer} and its accompanying materials would not have been possible without the dedicated efforts of the collaborators who worked tirelessly throughout Academic Year 2024. I would like to extend my deepest gratitude to undergraduate Instructional Assistants Kaylee Johnson (lead), Advaith (Adi) Balaji, Madeline (Maddy) Bezzina, Justin Boverhof, Anran (Annie) Li, Elaina Mann, Reina Mezher, and Maxwell (Max) West, whose invaluable insights and feedback were instrumental in shaping the homework sets, projects, and the initial draft of the textbook. I would also like to acknowledge the many contributions of Graduate Student Instructor Maxwell Gonzalez, whose reading of the textbook and solving all of the homework and projects have further refined the course. These young learners' diligence, creativity, commitment to testing and improving the materials, and supporting the students in the pilot semester have been essential to the successful offering of \textsc{Calculus for the Modern Engineer}. Finally, a special thank you to the group of students who were brave enough to take the pilot offering of ROB 201 \textsc{Calculus for the Modern Engineer} in Fall 2024. Your feedback and enthusiasm have been invaluable in shaping the future of this course.

\bibliographystyle{alpha}
\bibliography{ASEE2025}

\appendix
\setcounter{section}{0}
\renewcommand{\thesection}{Appendix \Alph{section}}:
\section{How the Course Material was Selected}
\label{sec:appendix}

\input{Sections/09ContentRationnale}

\section{Overview of ROB 101 - Computational Linear Algebra}
\label{sec:ROB101}
\input{Sections/10ROB101ComputationalLinearAlegbra}

\section{ROB 201 Midterm Teaching Evaluations}
\label{sec:ROB201midtermEvaluations}

\includegraphics[width=1.0\textwidth, page=1, trim=0 0 0 0, clip]{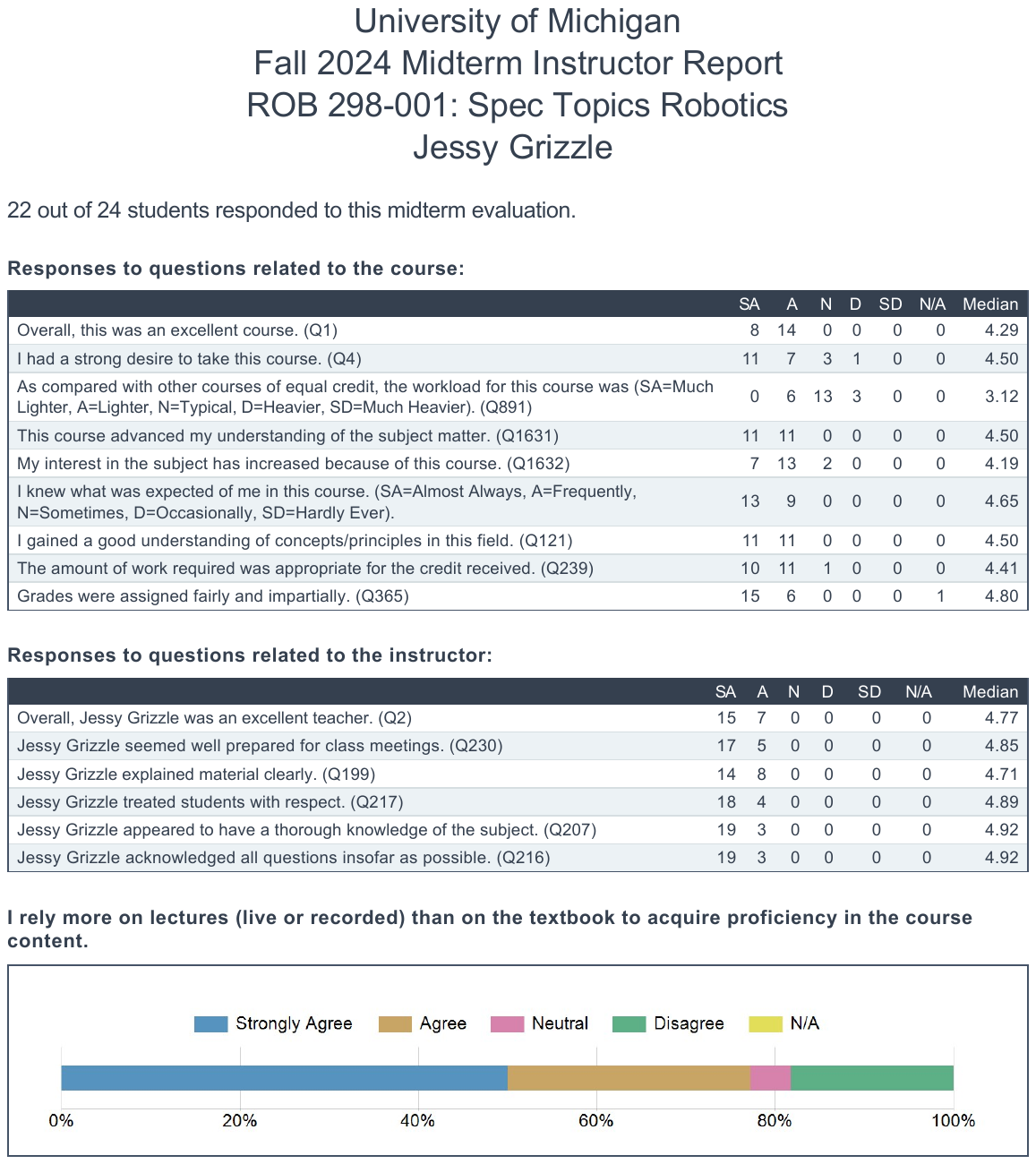}

\includepdf[pages=2-, width=1.0\textwidth, pagecommand={}]{PDFs/FA2024midterm.pdf}

\section{ROB 201 End-of-Term Teaching Evaluations}
\label{sec:ROB201endtermEvaluations}

\includegraphics[width=1.0\textwidth, page=1, trim=0 0 0 0, clip]{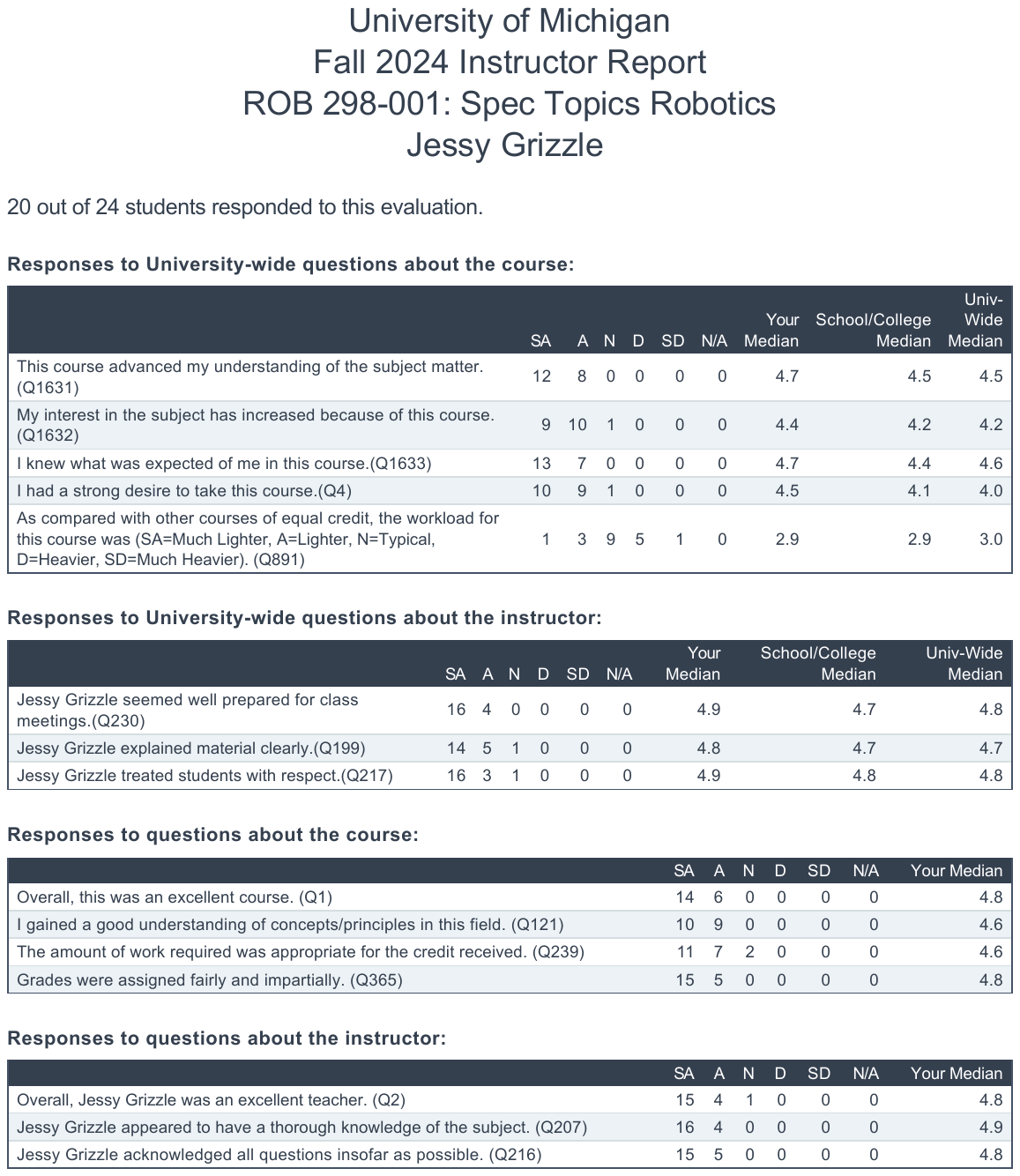}



\section*{Supplementary material}

# Appendix C   ROB 201 Midterm Teaching Evaluations

## University of Michigan
## Fall 2024 Midterm Instructor Report
## ROB 298-001: Spec Topics Robotics
## Jessy Grizzle

22 out of 24 students responded to this midterm evaluation.

**Responses to questions related to the course:**

| | SA | A | N | D | SD | N/A | Median |
|---|---|---|---|---|---|---|---|
| Overall, this was an excellent course. (Q1) | 8 | 14 | 0 | 0 | 0 | 0 | 4.29 |
| I had a strong desire to take this course. (Q4) | 11 | 7 | 3 | 1 | 0 | 0 | 4.50 |
| As compared with other courses of equal credit, the workload for this course was (SA=Much Lighter, A=Lighter, N=Typical, D=Heavier, SD=Much Heavier). (Q891) | 0 | 6 | 13 | 3 | 0 | 0 | 3.12 |
| This course advanced my understanding of the subject matter. (Q1631) | 11 | 11 | 0 | 0 | 0 | 0 | 4.50 |
| My interest in the subject has increased because of this course. (Q1632) | 7 | 13 | 2 | 0 | 0 | 0 | 4.19 |
| I knew what was expected of me in this course. (SA=Almost Always, A=Frequently, N=Sometimes, D=Occasionally, SD=Hardly Ever). | 13 | 9 | 0 | 0 | 0 | 0 | 4.65 |
| I gained a good understanding of concepts/principles in this field. (Q121) | 11 | 11 | 0 | 0 | 0 | 0 | 4.50 |
| The amount of work required was appropriate for the credit received. (Q239) | 10 | 11 | 1 | 0 | 0 | 0 | 4.41 |
| Grades were assigned fairly and impartially. (Q365) | 15 | 6 | 0 | 0 | 0 | 1 | 4.80 |

**Responses to questions related to the instructor:**

| | SA | A | N | D | SD | N/A | Median |
|---|---|---|---|---|---|---|---|
| Overall, Jessy Grizzle was an excellent teacher. (Q2) | 15 | 7 | 0 | 0 | 0 | 0 | 4.77 |
| Jessy Grizzle seemed well prepared for class meetings. (Q230) | 17 | 5 | 0 | 0 | 0 | 0 | 4.85 |
| Jessy Grizzle explained material clearly. (Q199) | 14 | 8 | 0 | 0 | 0 | 0 | 4.71 |
| Jessy Grizzle treated students with respect. (Q217) | 18 | 4 | 0 | 0 | 0 | 0 | 4.89 |
| Jessy Grizzle appeared to have a thorough knowledge of the subject. (Q207) | 19 | 3 | 0 | 0 | 0 | 0 | 4.92 |
| Jessy Grizzle acknowledged all questions insofar as possible. (Q216) | 19 | 3 | 0 | 0 | 0 | 0 | 4.92 |

**I rely more on lectures (live or recorded) than on the textbook to acquire proficiency in the course content.**

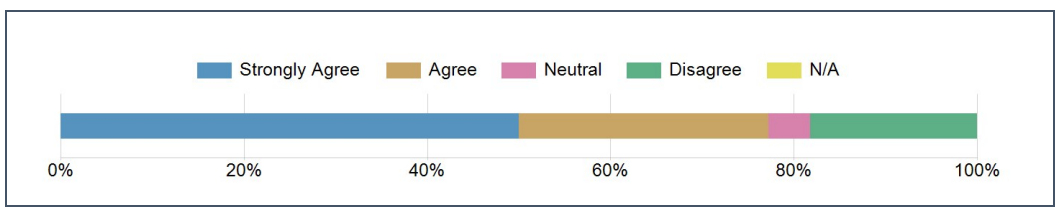

**The Julia programming exercises in this course have successfully helped bring alive the subject of calculus, Making it more relatable and practical, particularly in the context of engineering applications.**

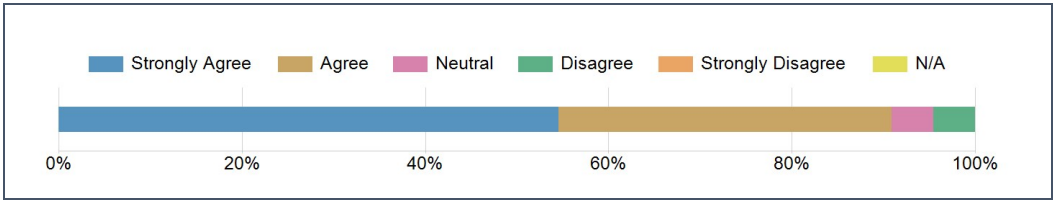

**The written homework assignments in this course have effectively helped me practice the mechanics of standard calculus operations, such as computing limits, finding areas, and differentiating functions.**

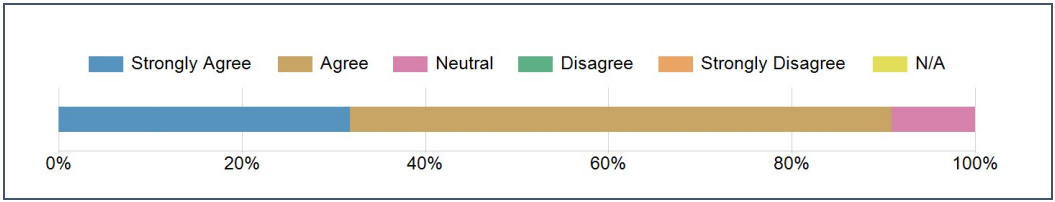

## Written Comments

**Which aspects of this course were most valuable? (Q908)**

| Comments |
|---|
| Learning calculus in a better way |
| Project and computation based approach to learning mathematics |
| I found the projects so far most valuable. |
| Julia homeworks |
| Dr. Grizzle is the most valuable aspect of the course. It's not often that you're being taught by the textbook author. He does a great job of explaining concepts, is very friendly, and I enjoy taking his classes. I do not think that this accelerated of a math course would be nearly as successful if taught by a staff member who didn't have a hand in creating this specific curriculum and textbook. |
| The most valuable aspects of the course in my opinion is the extremely fundamental approach that is taking in learning everything. We're never just told to "memorize x formula", but are shown the steps needed to get to said formula, giving a deeper understanding of various aspects of calculus. |
| I find it valuable how Professor Grizzle focused not just on teaching calculus but on demonstrating its real–world applications in robotics. |
| Seeing real world examples of how calculus is being used through the projects and programming homework |
| N/A |
| Learning actual mathematical applications of the math that we are learning and how it can be applied to robotics. |
| The multiple formats of the homework (written & Julia) was really great for being able to understand the math while also seeing the real applications of the math in code. The projects were also helpful for seeing these applications and made the course much more interesting and fun. |

**Which aspects of this course were least valuable? (Q909)**

| Comments |
|---|
| Nothing |
| Did not find anything least valuable. |
| Written Homeworks |
| I think the coding aspect of this course is the least valuable because while it might be the preferred way that engineers do math in industry, there is too much starter code for us to really learn any proper coding techniques. I had similar complaints with ROB 101. I was able to complete the assignments I had in Julia successfully, but walking away from that course I can confidently say I can't replicate any of the code on my own. |
| I would consider the less valuable aspects of the course to include the massive amount of content covered by the course. While it is a general overview, and combines topics from many math classes, the heavy load forces a rapid class pace, which can make keeping up difficult. |
| Nothing comes to mind> |
| Written homework currently feels like busywork |
| N/A |
| None |
| I found all of it very helpful and valuable. There are so many resources to take advantage of, but I personally did not attend office hours, just because I never had any questions that couldn't be solved through the lectures, textbook, or Piazza. |

**Please comment on the quality of the course as a whole. (Q911)**

| Comments |
| --- |
| I love it |
| The course is set up to provide enough practice to students through written homeworks projects and code homeworks. |
| Great |
| I think this course is really good for a starter course. I have no notes on the actual content of the class. I do think that it would benefit from a lab section specifically focused on Julia so that the coding section is more meaningful. This would allow less starter code, which encourages more students to learn coding techniques better. |
| I believe the support offered by the professor and the TAs supports the rapid pace of the course, and the passion shown by everyone involved contributes, from the custom written textbook to the piazza posts. |
| The course has be great –– I've learned a lot throughout the pass couple of weeks and I'm looking forward to the content we still have to learn in the future. |
| Very good! I think it's much more interesting and is more involved than a regular math class |
| Well planned and interesting |
| I think that the quality of the course as a whole has been very good so far for being a pilot semester, only minor issues have arisen and the content has so far made sense and the assignments we have been asked to do haven't been absurd or anything. |
| This course was very interesting and fun as well. I have really enjoyed this course as a way to learn calc in an engineering and robotics context instead of through the math department (sorry not sorry haha). |

**How can Jessy Grizzle improve the teaching of this course? (Q901)**

| Comments |
| --- |
| Hes the best |
| Have more examples in class that are similar to homeworks. |
| Offering more sections at different times |
| He is an excellent teacher. He has good pacing, explains high level concepts in low level manners making the class very approachable. And he is a personable guy so going to class is fun. |
| Not really sure. |
| I think Professor Grizzle is an amazing teacher and I really enjoy his teaching style. I don't think he has to change much to improve his teaching. |
| GRIZZLE IS THE BEST |
| Maybe go more in depth on the textbook than how it currently is for lectures |
| N/A |
| Occasionally in lecture, we cover topics too fast for me to write things down. I'd appreciate having a bit more time to catch up. |
| Explain some of the proofs a little better for the people who haven't really had any real introduction or work done with proofs before other than in this class. |
| If only more people attended lecture, but that's not his fault, the course is great! |

**Explain how Calculus for the Modern Engineer is meeting and/or not meeting your expectations as a new perspective on the subject matter:**

| Comments |
|---|
| It is meeting it I love the class |
| It is meeting expectations as I am able to learn applications of material to real world. |
| Met my expectations |
| This will just be me reiterating that for a class of this caliber that can count for up to 2 math courses, I feel that the practical application with Julia should be a bigger focus. We should be given data and minimal starter code. |
| I appreciate getting a review on a significant amount of calculus material, and being introduced to new concepts along the way. |
| This class is meeting my expectations because it's not just focused on purely theoretical calculus, it focuses on the practical applications in engineering. |
| The class is much more interesting than I thought it would be for a math class |
| It is meeting my expectations as it is an informative course covering the expected topics |
| This class is definitely meeting my expectations when it comes to giving a new perspective on calculus in general and how calculus is going to be used by engineers rather than a rigorous and unnecessarily busy math class with no applications or insight into how math is really being used by engineers. |
| As someone who has taken (and dropped) Math 216 in the past, I can say this course does a much better job at explaining the course material in a graspable and interesting way. I feel 10x more confident with even the material from Math 115 and 116 after seeing it in this course through the lectures and homeworks. |

**Discuss how the textbook's level of presentation is meeting or not your needs now, as well as its prospects as a long-term reference on the subject matter in terms of theoretical fundamentals and practical implementation of Calculus.**

| Comments |
|---|
| It is a well written textbook |
| The textbook is simple and encodes a lot of important information in color coded format making it easy for student to extract most important information. |
| Met my expectations |
| The textbook is aimed around the Julia programming language, and while this is helpful, I think including an appendix with python or C++ equivalents of the functions would make it a more useful book to turn back to considering the prevalence of these two languages in the industry. |
| My problem with the textbook is that it's a lot of (valuable) content, but it is difficult to peruse and use as a quick reference due to the large amount of content. The explanations are clear, and the links and graphs are useful, it just can be tricky to find specific topics. |
| The textbook explains the concepts well, but I wish it provided more thorough explanations of the code. A deeper breakdown of the logic, structure, and implementation would make it easier to understand how the code connects with the math concepts we are learning. |
| I think it does a very good job showing how different things that are learned in calculus are used in the real world, I have taken calc 1 and 2 in the past but I wouldn't have thought of how to use integration like how it's being used in the programming homeworks or projects. |
| The textbook is very clear and allows helpful practice problems |
| I think the textbook is very good, every time I read it I get a very good understanding of the content, and I really appreciate all of the implementation of functions and explanation of how to make some different functions or mathematical computations in Julia. |
| The textbook is such a great resource not only for derivations and equations, but for the applications and code snippets. If I ever have trouble with these calc concepts in the future, I am confident I could reference my ROB 201 textbook and be refreshed no problem. |

# Appendix D  ROB 201 End-of-Term Teaching Evaluations

University of Michigan
Fall 2024 Instructor Report
ROB 298-001: Spec Topics Robotics
Jessy Grizzle

20 out of 24 students responded to this evaluation.

**Responses to University-wide questions about the course:**

|  | SA | A | N | D | SD | N/A | Your Median | School/College Median | Univ-Wide Median |
|---|---|---|---|---|---|---|---|---|---|
| This course advanced my understanding of the subject matter. (Q1631) | 12 | 8 | 0 | 0 | 0 | 0 | 4.7 | 4.5 | 4.5 |
| My interest in the subject has increased because of this course. (Q1632) | 9 | 10 | 1 | 0 | 0 | 0 | 4.4 | 4.2 | 4.2 |
| I knew what was expected of me in this course.(Q1633) | 13 | 7 | 0 | 0 | 0 | 0 | 4.7 | 4.4 | 4.6 |
| I had a strong desire to take this course.(Q4) | 10 | 9 | 1 | 0 | 0 | 0 | 4.5 | 4.1 | 4.0 |
| As compared with other courses of equal credit, the workload for this course was (SA=Much Lighter, A=Lighter, N=Typical, D=Heavier, SD=Much Heavier). (Q891) | 1 | 3 | 9 | 5 | 1 | 0 | 2.9 | 2.9 | 3.0 |

**Responses to University-wide questions about the instructor:**

|  | SA | A | N | D | SD | N/A | Your Median | School/College Median | Univ-Wide Median |
|---|---|---|---|---|---|---|---|---|---|
| Jessy Grizzle seemed well prepared for class meetings.(Q230) | 16 | 4 | 0 | 0 | 0 | 0 | 4.9 | 4.7 | 4.8 |
| Jessy Grizzle explained material clearly.(Q199) | 14 | 5 | 1 | 0 | 0 | 0 | 4.8 | 4.7 | 4.7 |
| Jessy Grizzle treated students with respect.(Q217) | 16 | 3 | 1 | 0 | 0 | 0 | 4.9 | 4.8 | 4.8 |

**Responses to questions about the course:**

|  | SA | A | N | D | SD | N/A | Your Median |
|---|---|---|---|---|---|---|---|
| Overall, this was an excellent course. (Q1) | 14 | 6 | 0 | 0 | 0 | 0 | 4.8 |
| I gained a good understanding of concepts/principles in this field. (Q121) | 10 | 9 | 0 | 0 | 0 | 0 | 4.6 |
| The amount of work required was appropriate for the credit received. (Q239) | 11 | 7 | 2 | 0 | 0 | 0 | 4.6 |
| Grades were assigned fairly and impartially. (Q365) | 15 | 5 | 0 | 0 | 0 | 0 | 4.8 |

**Responses to questions about the instructor:**

|  | SA | A | N | D | SD | N/A | Your Median |
|---|---|---|---|---|---|---|---|
| Overall, Jessy Grizzle was an excellent teacher. (Q2) | 15 | 4 | 1 | 0 | 0 | 0 | 4.8 |
| Jessy Grizzle appeared to have a thorough knowledge of the subject. (Q207) | 16 | 4 | 0 | 0 | 0 | 0 | 4.9 |
| Jessy Grizzle acknowledged all questions insofar as possible. (Q216) | 15 | 5 | 0 | 0 | 0 | 0 | 4.8 |

**My favorite part of the course was the mathematical rigor it provided to help me learn the fundamentals of Calculus.**

(custom question added by the instructor)

| SA | A | N | D | SD | N/A | Your Median |
|---|---|---|---|---|---|---|
| 5 | 11 | 1 | 2 | 0 | 0 | 4.1 |

**My favorite part of the course was the applications it explored to help me understand the utility and power of Calculus in real-life engineering.**

(custom question added by the instructor)

| SA | A | N | D | SD | N/A | Your Median |
|---|---|---|---|---|---|---|
| 13 | 7 | 0 | 0 | 0 | 0 | 4.7 |

**You were glad to have been given a free 1-year license to Wolfram Alpha Pro because it will help you to do better in your engineering courses that use Calculus.**

(custom question added by the instructor)

| SA | A | N | SD | D | N/A | Your Median |
|---|---|---|---|---|---|---|
| 16 | 3 | 1 | 0 | 0 | 0 | 4.9 |

The medians are calculated from Fall 2024 data. University-wide medians are based on all UM classes in which an item was used. The school/college medians in this report are based on classes that are lower division with enrollment of 16 to 74 in College of Engineering.

## Written Comments

**Which aspects of this course were most valuable? (Q908)**

| Comments |
|---|
| I feel like the most valuable parts of this course were the programming assignments. As a neurodivergent student, I find it really difficult to complete tasks without a form of immediate feedback telling me I'm on the right track, and the friendly checks in the programming assignments helped a lot with this. |
| The projects, as it let me see how to actually use math in relation to problems |
| Actually applying calculus through various projects, and the ease of asynchronous learning. |
| I really appreciated the programming aspect including the programming homeworks and projects. I felt like they did a good job of providing context where the math concepts would be applied and also will help me with future classes and the real world. |
| The multiple formats for the homework helped me to understand the material and know how to apply it to real work scenarios. |
| The most valuable aspects of this course were the depth and clarity of the material, as well as the way it bridged theoretical concepts with practical applications. |
| PROF. GRIZZLE!!!!!! |
| The experience with software tools to evaluate and solve math problems. |
| Learning the applications of how calculus and the concepts we are learning can be applied and how we might be able to use them in the future in either other classes or in our future professional fields. |
| Julia Homework |
| The most valuable aspects were the addition of written questions to test mathematical concepts combined with coding questions for practical use. |
| Julia homeworks |

**Which aspects of this course were least valuable? (Q909)**

| Comments |
|---|
| written homeworks |
| I found it all valuable |
| The least valuable aspect of this course was the scheduling, as the early class time made it challenging for many students to attend and fully engage. |
| Some of the julia homework assignments felt kind of unnecessary, and I would have rather learned a little more by hand rather than had to move on from doing written homework over to Julia homework. |
| Written HW |
| There were not any. |
| None |

**Please comment on the quality of the course as a whole. (Q911)**

| Comments |
|---|
| It was good! There were some hiccups since it was the first year for the class, however they were ironed out quickly. |
| The quality was very high. I don't know any other courses where the person teaching the course wrote the textbook and is as knowledgeable about the subject as this one. |
| The course is shockingly high quality for a course that has never happened before. I was impressed with the breadth of resources available, and I shockingly actually liked the textbook. |
| This has been an excellent course and I would highly recommend it to anyone interested/required to take differential equations |
| The course as a whole is highly valuable and well–structured. It offered me thorough understanding of the material and its real–world applications. The content is challenging yet rewarding, pushing me to deepen my problem–solving skills and critical thinking when it comes to calculus. |
| I think the class so much merit it should replace Math 216 in our curriculum. I think if the credits the class counted for were adjusted to Math 215, and Math 216, students coming in from the base calc sequence would still learn a lot. Given how much other classes rely on skills learned in Math 115 and 116 (especially in exams) it would benefit students more to be taught the real engineer's way to solve 115 and 116 style problems computationally. Adjusting the course like this would reduce the amount of time needed to cover the 115 and 116 material, and allow for expanding the 215 and 216 material. I think those changes would make Rob 201 an excellent foundation for the much more difficult math required in 300s and 400s in the ROB department. |
| I think the course as a whole was great, it helped me learn a lot, and it helped to enlighten me on a lot of the important things I could use later or things that don't get taught in a normal math course. |
| It was great |
| It was a great course which taught in a way that will be useful for many robotics majors. |
| Absolutely loved it. Thank you! |
| Excellent course |

**How can Jessy Grizzle improve the teaching of this course? (Q901)**

| Comments |
|---|
| I thought pacing of lecture was good, no changes needed |
| If possible, I think shorter lectures would be nice so it's easier to pay attention. I would rather have two 1.5 hour lectures and one 1 hour discussion a week. That being said, the course was really a very good one. |
| I think this course would just be better with a bit more participation from students, and a larger attendence. |
| Professor Grizzle has done an excellent job teaching this course, particularly by providing a remote option, which added flexibility for anyone taking the class. However, the early class timing seemed to impact student engagement over time. |
| Professor Grizzle is an excellent lecturer, the remote options and lecture recordings make the course accessible to a variety of learners. |
| I don't think there is much that can be done to improve his teaching of the course, I loved the enthusiasm that he brought to the classroom, and also the extensive knowledge of robotics and math that helped me learn. The only critique I would have for the teaching aspect, would be to explain concepts further than just a proof, to me at least, proofs can be quite confusing, so I think it would be helpful to either teach proofs better at the beginning of the course, or supplement more example problems that cover everything in order to teach concepts as I feel like, I personally learn better off of an example rather than a proof. |
| Nothing |
| Nothing. |
| By continuing to bring his great energy to the course. |

**If you like the course, give a pitch you would make to entice other students to also take it. If you dislike the course, give a pitch to dissuade other students from taking it.**

| Comments |
|---|
| Class is pretty easy and gives a lot of helpful knowledge for actual usage of calculus |
| This is course is definitely difficult, but if you want to learn a large spread of calculus concepts in 1 semester, there's no better way to do it. The teacher and TAs are almost begging for students to use resources like office hours, and they're extremely responsive to any questions. Also, the class supports completely asynchronous learning, which as someone with a very busy schedule, I found it extremely helpful. |
| Calculus doesn't need to be a dreaded prerequisite before you take the "interesting classes". This class has not only made learning calculus concepts more engaging and enjoyable, but it has given me the tools to apply what I've learned to real engineering problems. |
| The only drawback I found in this course was its scheduling. A 10:30 AM class time might be challenging for students, and I believe shifting it to a later time in the day could significantly improve attendance and engagement. Professor Grizzle creates an engaging learning environment and adjusting the schedule would allow even more students to fully participate and benefit from the course. |
| If you've already taken math 115 and 116, 201 is going to teach you everything you actually need to know about Diff–eq to solve the problems the 300 levels will throw at you. It will also give you insight into robot controls. |
| If you are interesting in math, robotics, or the applications of math, this course is definitely for you. Personally I would take this course just to learn the applications of math and how these concepts relate to both the field of robotics, and the world arounds us. This class also teaches you about things that can be useful in other classes, or concepts that can help you be more successful as an engineer that aren't taught in a typical math class. Along with that you get to learn a programming language better, and can learn how a programming language can do a lot of the hard work of math for you, in order to make your life as an engineer easier. |
| It was very flexible |
| I would encourage others to take this course because it is not overwhelming and is able to teach more practical implementations of concepts we spend lots of time in calculus math class. |
| Have your calculus classes felt uninteresting, abstract, or made you feel like you will never be good at math? If so, take this course. |
| This class provides a very interesting and modern approach to Calculus, which I found more applicable to my real life/degree than any other math course I've taken at UMich. There is also opportunity for remote participation, which I took advantage of (and was very helpful when trying to balance a heavy course load). I strongly encourage taking this class as an alternative to MATH 116/216, as it is far more understandable and interesting. Professor Grizzle is great, explains concepts clearly, and tries to help students understand material as much as possible. The textbook is also well–written and engaging and helped me to understand material. |

**Computation was primarily used in HWs to show that the oft-perceived drudgery of Calculus could be avoided with the correct software tools, and in the Projects to show that realistic problems could be approached with knowledge from the course. If you were to re-design how computation were used in ROB 201, how would you go about it? Keep in mind that most students can only spend 4 hours per week doing HW for any given course.**

| Comments |
|---|
| It felt like the written homeworks were just busy work, maybe finding some way to make them feel more meaningful? |
| I would not change anything about the computation aspect super dramatically. I really liked how computation was used in ROB 201 |
| I think there is an opportunity to assign computation problems to the class as a way to prepare/participate in recitations. In a larger class problems could be assigned to groups of students to work out together in class, or to be taken home and tried independently, and then review the answers during recitation. This could be for credit or not, and would give students an opportunity to work through problems and then come to recitation or lecture ready with questions if they got stuck on something. |
| I believe the way the class is set up now is perfect. I was spending at least 4 hours on all my homeworks per week. |
| Continuing with my idea that future students will already have taken 115 and 116, I would probably make the introductory calculus homeworks based on difficult web–work problems and exam questions from previous years. I think 2–3 weeks of students solving their demon problems from previous semesters in minutes instead of hours would be a great way to introduce them to the software tools. |
| I personally like the computation portion of the course, I feel like it forces you to learn fundamental skills that are being taught in the class. If I were to suggest something, I would say to keep the questions about making a cheat sheet or such, as it helps make sure you go back and review the topics in the book. I personally don't have any complaints at the moment about the computation portion of the class. I do want to say that the recitations that were done in class did make the computation homework a lot less intimidating, since there was a part of class to teach you specifically about that topic since it was coming up in the homework. |
| I wouldnt |
| I would not redesign the hws or projects. |
| I liked how it was used right now. I wouldn't redesign it.f |
| I wouldn't redesign it, I liked this aspect of the class |



\end{document}

%% file: Defines/definesPackages.tex
\usepackage{yfonts}

\usepackage{paralist}

\usepackage{array}
\usepackage{scalerel}
\usepackage{setspace}
\usepackage{amsfonts}
\usepackage{mathtools} 
\usepackage{bm}
\usepackage{upgreek}

\usepackage{empheq}

\usepackage{tocloft}

\usepackage{float}
\usepackage{alltt}

\usepackage{placeins}

\usepackage{url}
\usepackage[utf8]{inputenc} 
\usepackage[T1]{fontenc}
\usepackage{lmodern}

\usepackage{pdfpages}
\usepackage{enumerate}
\usepackage[most]{tcolorbox}
\tcbuselibrary{theorems}

\usepackage{siunitx}

\usepackage{MnSymbol}

\usepackage{mathdots}
\usepackage{xfrac}

\usepackage{avant} 

\usepackage{subcaption}
\usepackage{graphicx} 
\graphicspath{{graphics/}}
\usepackage{amsthm}
\usepackage{algorithmicx}
\usepackage{mathrsfs}
\usepackage[ruled,vlined]{algorithm2e}
\usepackage{times}
\usepackage{multicol}
\usepackage[T1]{fontenc}
\usepackage{textcomp}
\usepackage{balance}
\usepackage{multirow}
\usepackage{booktabs}

\usepackage{makecell}
\usepackage{hhline} 
\usepackage[colorlinks]{hyperref}
\definecolor{BrickRed}{rgb}{0.7,0.15,0.15}
\hypersetup{colorlinks=true,pdfpagemode=UseNone,citecolor=black,linkcolor=black,urlcolor=BrickRed}

\definecolor{mylightblue}{rgb}{0.8,0.9,1}
\usepackage{comment}
\usepackage{stackengine}

\usepackage{xparse}
\usepackage{lipsum}
\usepackage{soul} 
\usepackage{times}
\definecolor{shadecolor}{rgb}{0.9,0.9,0.9}
\usepackage{enumitem}

\usepackage{mdframed}
\usepackage{framed}

\usepackage{listings}

\definecolor{codegreen}{rgb}{0,0.6,0}
\definecolor{codegray}{rgb}{0.5,0.5,0.5}
\definecolor{codepurple}{rgb}{0.58,0,0.82}
\definecolor{backcolour}{rgb}{0.95,0.95,0.92}

\lstdefinestyle{mystyle}{
    backgroundcolor=\color{backcolour},   
    commentstyle=\color{codegreen},
    keywordstyle=\color{magenta},
    numberstyle=\tiny\color{black},
    stringstyle=\color{black},
    basicstyle=\ttfamily\normal,
    breakatwhitespace=false,         
    breaklines=true,                 
    captionpos=b,                    
    keepspaces=true,                 
    numbers=left,                    
    numbersep=5pt,                  
    showspaces=false,                
    showstringspaces=false,
    showtabs=false,                  
    tabsize=2
}
\usepackage{mathtools}
\usepackage[english]{babel}


%% file: Defines/definesCommands.tex
\newcommand{\real}{\mathbb R}

\DeclareDocumentCommand{\zeros}{ O{} }{\textbf{0}_{#1}}

\usepackage{upgreek}

\definecolor{darkblue}{rgb}{0,0,0.5}
\definecolor{lightblue}{rgb}{0.88,1,1}
\definecolor{jwgbluegray}{rgb}{0.8, 0.95, 0.95} 

\definecolor{none}{cmyk}{0.0, 0.0, 0.0, 0.0}
\definecolor{brightblue}{RGB}{0, 0, 239} 

\newtcolorbox{emstatbox}[1][]{
  breakable,
  colback=jwgbluegray,
  colframe=brightblue, 
  boxsep=0pt,
  left=2pt,
  right=2pt,
  top=5pt,
  bottom=5pt,
  arc=0pt,
  outer arc=0pt,
  boxrule=2pt, 
  #1
}

\newcommand{\emstat}[1]{\begin{center}\begin{emstatbox}#1\end{emstatbox}\end{center}}

\definecolor{note}{rgb}{0.3,0.7,0.25}
\definecolor{rephase}{rgb}{0.15,0.7,0.15}
\definecolor{bag}{rgb}{0.6,0.6,0.2}

%% file: Sections/01Intro.tex
The calculus curriculum taught at most universities today remains rooted in a structure developed during the 1950s. Designed during the Sputnik era, this traditional sequence was intended to meet the educational needs of its time, emphasizing manual problem-solving techniques and abstract theoretical concepts. However, the landscape of engineering has evolved dramatically since then. Fields such as robotics, artificial intelligence, and data science now require not only a deep theoretical understanding but also the ability to engage confidently with computational tools and solve real-world problems. Despite these changes, the core structure of calculus courses---typically covering differential and integral calculus along with ordinary differential equations---remains largely unchanged, creating a widening gap between classroom learning and the dynamic problems modern engineers face.

To address this gap, the University of Michigan's Robotics Department launched an experimental course, \textsc{Calculus for the Modern Engineer}. This pilot course is part of an ongoing effort to explore how mathematics can be better aligned with the computational and practical needs of students entering fields like robotics and engineering. The goal of this article is to report on a new curriculum that integrates computation and real-world applications into calculus education while simultaneously making it a more engaging subject to learn.

A key challenge the course seeks to address is how calculus is traditionally taught. The focus on rote memorization and mechanical application of formulas often stifles curiosity, skill acquisition, and deep understanding. Students are rarely shown how calculus can be a creative, powerful tool for solving real-world engineering problems. Instead, they are burdened by repetitive manual computations that fail to connect with the complex, technology-driven tasks they will face in their careers. Moreover, outdated assessment methods, such as high-stakes exams, reinforce this disconnect, often producing median scores below 40\% and fostering frustration rather than confidence.

This experimental course aims to reframe calculus education to emphasize not just theoretical principles but their practical applications. Traditional timed exams have been replaced with three major projects that focus on solving real-world engineering challenges using the Julia programming language. This shift in assessment methods reflects the evolving needs of modern engineers, who increasingly rely on computational tools and iterative algorithms rather than closed-form solutions. The course allows students to engage with calculus concepts by directly applying them to projects, such as numerically integrating drone data, optimizing robotic movement, and designing feedback control systems.

By introducing computational tools early and focusing on real-world applications, the pilot course strives to rekindle students' enthusiasm for calculus. The course leverages Julia to help students visualize and manipulate complex mathematical models, fostering a deeper, more intuitive understanding of calculus concepts. Rather than being seen as a barrier, calculus becomes a dynamic enabler of engineering solutions.

This article reports on the progress and outcomes of this pilot course at the University of Michigan, detailing its structure, the rationale behind its development, and the practical experiences of both students and instructors. While the experiment is ongoing, early feedback suggests that this reimagined approach may offer a more relevant and engaging way to teach calculus to the engineers of tomorrow.

%% file: Sections/02ReformEfforts.tex
Here we give a snapshot of prior and current reform efforts.

\subsection{Calculus in Context}

\textit{Calculus in Context}, developed by the Five College Consortium\footnote{Smith College, Amherst College, Hampshire College, Mount Holyoke College, and the University of Massachusetts Amherst.} in 1995, restructured the traditional calculus sequence to make the subject more engaging, de-emphasize hand computations, and embrace computational methods \cite{callahan1995calculus}. Their innovative approach aimed to restore calculus as a dynamic tool for scientific inquiry, akin to how it was used by pioneers like Euler and Bernoulli. Key elements of their reform included:

\begin{enumerate}
    \item \textbf{Making Calculus Engaging Again}: The curriculum treated calculus as a powerful tool for exploring real-world problems, driven by scientific and mathematical questions that made the learning process more meaningful.
    
    \item \textbf{Replacing Hand Computations with Algorithms}: Embracing the potential of technology, the program prioritized numerical methods and algorithmic solutions over closed-form solutions, which are limited to special cases. This encouraged students to approach calculus iteratively through computational tools.
    
    \item \textbf{Deprecating Special-Case Topics}: The curriculum downplayed topics with narrow applications in favor of broadly applicable concepts such as differential equations and dynamical systems, helping students apply calculus to model real-world phenomena.
    
    \item \textbf{Collaborative and Experimental Learning}: A key pedagogical innovation was fostering a collaborative environment where students tackled messy, real-world problems, gaining an appreciation for the value of approximate solutions and the discovery process.
    
    \item \textbf{Shifts in Emphasis}: The curriculum focused on concepts, geometry, graphs, and numerical methods rather than techniques, formulas, and closed-form solutions. Students learned to appreciate when brute-force methods like Euler’s method were the appropriate choice over more elegant but less practical approaches.
\end{enumerate}

Despite these innovations, the curriculum has not been widely adopted. The University of Redlands \cite{beery_CiC_info} found challenges in the curriculum. Issues like unpolished materials, assumptions about students’ prerequisite knowledge, and articulation problems between Calculus II and III contributed to a partial return to traditional calculus curricula. Nevertheless, the faculty noted improved conceptual understanding, especially regarding limits and derivatives, and remained enthusiastic about the course’s potential for discovery-based learning.

\subsection{Calculus for Economists}

In the textbook, \texttt{Calculus for Economists} \cite{ahmed2024calculusTextbook}, published in September 2024, Ahmed integrates foundational topics from Calculus I, Calculus II, and Ordinary Differential Equations into a comprehensive curriculum tailored specifically for students in economics. The course associated with the textbook is structured into three primary modules \cite{ahmed2024calculusModuleI}, beginning with an introduction to calculus principles, followed by an in-depth exploration of differential calculus concepts and their practical applications. The final segment delves into advanced techniques of differential calculus, directly linking calculus concepts with key economic concepts, including price determination, consumer behavior, production theories, and fundamental macroeconomic principles. 

By bridging these fields, the textbook serves as a resource connecting calculus in economic modeling and analysis, thereby reinforcing the practical relevance of mathematical skills in economic decision-making and policy formulation. 
The curriculum is designed to enhance students' mathematical skills while also deepening their understanding of economic principles through applied calculus.

\subsection{Other Collaborative Efforts to Innovate Calculus Education}

\textit{SimCalc: Accelerating Access to the Mathematics of Change} represents an ambitious effort to democratize calculus by introducing core concepts of \textit{change and variation}---such as rates, accumulation, and limits---to students well before they encounter traditional calculus courses \cite{hegdeus2013simcalc}. Through SimCalc's “MathWorlds” software, students engage with dynamic visual representations of real-world motion scenarios (e.g., velocity and acceleration), fostering both qualitative and quantitative reasoning. This initiative aims to address disparities in mathematical education by providing interactive, technology-based learning environments that make abstract calculus concepts accessible to a broader range of students.

While SimCalc has developed powerful tools to support learning, including dynamic graphs and visualizations, the project is still in the process of refining its curriculum for broader classroom use. Its goal is to not only transform how calculus is taught but also to increase participation and engagement among underrepresented groups, ultimately integrating technology, pedagogy, and curriculum reform in mathematics education.

\textit{Learning Lab Calculus Report} by the California Learning Lab emphasizes the critical gatekeeping role that calculus plays in determining STEM persistence, particularly for underrepresented student populations \cite{calearninglab2021calculus}. The report calls for active learning methods, contextualized teaching, and pedagogical shifts that address systemic disparities in K–12 preparation and promote inclusivity in college-level mathematics. Notable recommendations include integrating mentorship programs, fostering collaborative learning environments, and implementing adaptive teaching practices to better support a diverse range of learners.

\textit{Project NExT} (New Experiences in Teaching), sponsored by the Mathematical Association of America (MAA) \cite{ProjectNExT}, is a comprehensive professional development program for new mathematics faculty. While not exclusively focused on calculus, it seeks to significantly impact mathematics education across various levels by promoting teaching innovations, active learning techniques, and learning technology integration in mathematics instruction. Project NExT addresses curriculum development, emphasizing student engagement and relevance, particularly in core courses like calculus. It also tackles critically important issues such as equity in mathematics education, supporting diverse student populations, and balancing teaching responsibilities with research. By creating a supportive network of early-career mathematicians, Project NExT facilitates the sharing of best practices and resources across institutions. 

It goes almost without saying that much remains to be done to meet these laudable goals \cite{bressoud2019calculus}. \textsc{Calculus for the Modern Engineer} hopes to be a significant step forward. 

\subsection{Efforts Toward Interactive Calculus Texts or the Inclusion of Programming}

Interactive textbooks that combine explanations with live programming cells provide a potentially transformative approach to teaching calculus, enabling students to actively engage with the material. Noteworthy initiatives include:

\begin{itemize}
    \item \texttt{Active Calculus} is a free, open-access textbook that breaks away from passive learning by encouraging students to explore calculus through computational experiments and real-world applications \cite{activecalculus}. 
    \item \texttt{Calculus With Julia} offers computational tools for interactive learning of calculus principles and has been used in \texttt{MTH 229: Calculus Computer Laboratory} \cite{verzani2014mth229, verzani2024calculuswithjulia}.
\end{itemize}

These efforts build upon the pioneering work of \textit{Calculus in Context}, discussed earlier. Despite their significant contributions, they adhere to the traditional sequence of calculus topics, which may limit their ability to fully revolutionize the curriculum. True transformation requires rethinking not just how calculus is taught, but also which topics are emphasized and in what order they are introduced.

The incorporation of computational tools into mathematical education across various disciplines seeks enhanced understanding and application of complex concepts. For instance, Brokate et al. (2019) highlight the integration of MATLAB in teaching calculus to engineers, providing students with real-world problem-solving skills \cite{brokate2019calculus}. Similarly, Adhikari et al. (2015) emphasize the importance of computational thinking in data science, offering foundational skills in inferential thinking that are crucial for modern data-driven environments \cite{adhikari2015computational}. Furthermore, Lipsman and Rosenberg (2017) and Brown et al. (1991) have also explored the use of MATLAB and Mathematica to teach multivariable calculus and reinforce its applications in geometry and physics \cite{lipsman2017multivariable, brown1991calculus}.

\subsection{Comparison and Contrast with ``Calculus for the Modern Engineer''}

\textsc{Calculus for the Modern Engineer} aligns with the educational philosophies seen in \textit{Calculus in Context} \cite{callahan1995calculus} and \textit{Calculus for Economists} \cite{ahmed2024calculusTextbook}, both of which also compress the traditional calculus sequence into a comprehensive curriculum infused with real-world examples. This alignment underscores a shared emphasis on enhancing the practical relevance of calculus, thereby maintaining the discipline's significance in contemporary education.

Distinctively, \textsc{Calculus for the Modern Engineer} not only embraces application integration but also reimagines the curriculum by reordering topics into a sequence that progressively builds student skills in a manageable and coherent manner. This methodical restructuring allows for deeper comprehension and sustained knowledge retention, aiding students in their transition to more complex mathematical and engineering concepts. The course uniquely prepares students for advanced topics in robotics by emphasizing a solid theoretical foundation alongside practical problem-solving skills; most of the results in the textbook are also proved in optional read sections. This overall approach strives to equip students with a robust mathematical toolkit tailored to meet the challenges of modern technological and engineering landscapes.

%% file: Sections/03TextbookContent.tex
The textbook accompanying \textsc{Calculus for the Modern Engineer} is a comprehensive resource, covering everything from basic pre-calculus principles to advanced topics like ODEs, Laplace transforms, and proportional-derivative feedback controllers. The book places equal weight on the conceptual and computational sides of calculus. Concepts are first grounded in intuition, followed by epsilon-delta arguments and limits as needed, and are always closely followed by computations in Julia that illustrate the ideas and implement the associated formulas.

The material begins with a review of notation, functions, and algebraic manipulation before diving into the heart of calculus: definite integration. After covering the Riemann integral, the text explores differentiation, symbolic computation, and practical applications such as optimization with constraints. Students learn to use tools like SymPy and Wolfram Alpha Pro to perform symbolic manipulations, and they apply these concepts to robotics case studies involving mobile robots. The textbook concludes with chapters on improper integrals, ordinary differential equations, and Laplace transforms, the latter framed through the lens of feedback control and modern robotics applications.

\begin{enumerate}

\item {\bf Pre-calculus: Notation, Functions, and Various Algebraic Facts}  

\textbf{Notes:} $\bullet$ 2-hour lecture. $\bullet$ Students are expected to review this material mostly on their own. $\bullet$ 20-question quiz to ensure mastery. $\bullet$ The Approximation Principle is highlighted in HW01.  
\textbf{Learning Objectives:} $\bullet$ Recognize the utility of mathematical notation for precision and expressivity. $\bullet$ Understand and apply the Bisection Algorithm as an example of the Approximation Principle. $\bullet$ Reaffirm understanding of fundamental concepts such as functions, domains, ranges, and inverses.  
\textbf{Content:} Introduction, Got Calculus Dread? Additional Resources, Notation or the Language of Mathematics, The Approximation Principle, Algebraic Manipulation and Inequalities, Functions (Domains, Ranges, Inverses, and Compositions), Trigonometric and Inverse Trigonometric Functions, Powers and Roots, Exponentials and Logarithms, Euler’s Formula, Hyperbolic Trigonometric Functions, Summing Symbol, Binomial Theorem, Special Functions, Shifting and Scaling, and Feedback Control. (Optional Reads: Binomial Theorem meets Euler and Proofs Associated with the Chapter.)  

\item {\bf Calculus Foundations: Proofs, Finite Sums, Limits at Infinity, and Geometric Sums}  

\textbf{Notes:} $\bullet$ 4 hours lecture + recitation. $\bullet$ Written and Julia HWs.  
\textbf{Learning Objectives:} $\bullet$ Understand the significance of mathematical proofs. $\bullet$ Master proof by induction and apply it to sums of powers of integers. $\bullet$ Comprehend limits and their initial uses in calculus, with emphasis on limits at infinity for rational functions and exponentials.  
\textbf{Content:} Introduction, Mathematical Proofs, Countable Sets, Proofs by Induction, Finite Sums, Limits at Infinity, Geometric Sums. (Optional Reads: Euler's Number, Proofs Associated with the Chapter.)  

\item {\bf Definite Integration as the Signed Area Under a Curve}  

\textbf{Notes:} $\bullet$ 5 hours lecture + recitation. $\bullet$ Written and Julia HWs.  
\textbf{Learning Objectives:} $\bullet$ Define and explain definite integration \`a la Riemann-Darboux. $\bullet$ Execute computations using methods like Trapezoidal Rule and Simpson's Rule. $\bullet$ Recognize applications of definite integrals in engineering.  
\textbf{Content:} The Riemann Integral, Properties of the Riemann Integral, Numerical Methods for Approximating Integrals, Applications of Definite Integral. (Optional Reads: Proofs Associated with the Chapter.)  

\item {\bf Properties of Functions: Left and Right Limits, Types of Continuity, Boundedness, and Generalizations of Max and Min}  

\textbf{Notes:} $\bullet$ 5 hours lecture + recitation. $\bullet$ Written and Julia HWs.  
\textbf{Learning Objectives:} $\bullet$ Analyze function behavior using one-sided limits. $\bullet$ Understand the nature of function continuity. $\bullet$ When can limits be taken inside a function. $\bullet$ Explore generalizations of maximum and minimum values.  
\textbf{Content:} Limits from the Left and Right, Continuity, Maximum and Minimum Values, The Squeeze Theorem, Piecewise Continuity, Boundedness, Generalizations of Max/Min. (Optional Reads: Intermediate Value Theorem, Proofs Associated with the Chapter.)  

\item {\bf Differentiation}  

\textbf{Notes:} $\bullet$ 5 hours lecture + recitation. $\bullet$ Written and Julia HWs.  
\textbf{Learning Objectives:} $\bullet$ Understand single-variable and partial derivatives. $\bullet$ Apply differentiation rules and software tools. $\bullet$ Explore real problems where derivatives are used in engineering.  
\textbf{Content:} Derivative as Local Slope and as Linear Approximation, Differentiation Rules, Software Tools for Computing Derivatives, Use Cases of Derivatives, Partial Derivatives, Jacobians, Gradients, Hessians, Total Derivative. (Optional Reads: Proofs Associated with the Chapter.)  

\item {\bf Engineering Applications of the Derivative}  

\textbf{Notes:} $\bullet$ 6 hours lecture + recitation. $\bullet$ Written and Julia HWs.  
\textbf{Learning Objectives:} $\bullet$ Develop strategies for solving optimization problems. $\bullet$ Apply Lagrange’s equations to dynamics.  
\textbf{Content:} Path Length, Root Finding, Unconstrained Minimization, Gradient Descent with and without Equality Constraints, Lagrange Multipliers, Dynamics à la Lagrange. (Optional Reads: More on Lagrangian Dynamics, Proofs Associated with the Chapter.)  

\item {\bf Antiderivatives and the Fundamental Theorems of Calculus}  

\textbf{Notes:} $\bullet$ 5 hours lecture + recitation. $\bullet$ Written and Julia HWs.  
\textbf{Learning Objectives:} $\bullet$ Understand antiderivatives and their relationship to definite integrals. $\bullet$ Apply key techniques for finding antiderivatives.  
\textbf{Content:} Introduction, Fundamental Theorems of Calculus, Antiderivatives, Techniques for Finding Antiderivatives, Software Tools for Antiderivatives. (Optional Reads: Proofs Associated with the Chapter.)  

\item {\bf Improper Integrals}  

\textbf{Notes:} $\bullet$ 2 hours lecture + recitation. $\bullet$ Written and Julia HWs.  
\textbf{Learning Objectives:} $\bullet$ Define and understand improper integrals. $\bullet$ Explore practical applications in probability theory.  
\textbf{Content:} Type-I Improper Integrals, Type-II Improper Integrals, Improper Integrals in Probability. (Optional Reads: Proofs Associated with the Chapter.)  

\item {\bf Ordinary Differential Equations (ODEs)}  

\textbf{Notes:} $\bullet$ 9 hours lecture + recitation. $\bullet$ Written and Julia HWs.  
\textbf{Learning Objectives:} $\bullet$ Learn analytical and numerical methods for solving ODEs. $\bullet$ Apply ODEs to model dynamic systems in robotics.  
\textbf{Content:} Introduction to ODEs, Higher-Order ODEs, Vector ODEs, Linear Systems of ODEs, Numerical Solutions, Matrix Exponential, Eigenvalues and Eigenvectors. (Optional Reads: Resonance in ODEs, Proofs Associated with the Chapter.)  

\item {\bf Laplace Transforms through the Lens of Feedback Control}  

\textbf{Notes:} $\bullet$ 9 hours lecture + recitation. $\bullet$ Written and Julia HWs.  
\textbf{Learning Objectives:} $\bullet$ Understand Laplace transforms and their application in solving ODEs. $\bullet$ Derive and apply transfer functions in feedback control systems.  
\textbf{Content:} Laplace Transform, Solving ODEs Using Laplace, Transfer Functions, Poles and Zeros, Feedback Systems, Transient Response of Systems, Proportional and Proportional-Derivative Feedback, Pre-compensators. (Optional Reads: Impulse Function, Proofs Associated with the Chapter.)

\end{enumerate}

%% file: Sections/04Projects.tex
The structure of \textsc{Calculus for the Modern Engineer} is supported by three major Julia-based projects, each designed to emphasize core skills while progressively building students' computational and theoretical knowledge. These projects are carefully aligned with the needs of robotics and engineering, providing students with hands-on experience that bridges mathematical concepts with real-world applications. Each project introduces students to key tools and methodologies used in modern engineering, ensuring they are equipped to tackle complex problems using both analytical and computational approaches. There are eight homework assignments where more basic skills are acquired and honed. Each assignment has two parts, a standard written component mirroring traditional expectations in hand computations, and a jupyter notebook written in Julia. All of this material will be open-sourced~\cite{grizzle2025rob201}.

\subsection{Numerical Integration: Estimating Speed and Position with IMU Data}

The first project focuses on numerical integration, a fundamental concept in calculus and robotics. Students are tasked with using linear acceleration data collected from the IMU of a drone to estimate the drone’s velocity and position. Initially, they work with clean data, applying the Trapezoidal Rule to build a real-time estimate of speed and position:
$$ v(t) = v(t_0) + \int_{t_0}^t a(\tau) \, d \tau, $$
$$ p(t) = p(t_0) + \int_{t_0}^t v(\tau) \, d \tau. $$

As the project progresses, students encounter the real-world complication of acceleration bias, which affects the accuracy of the estimated speed and position. They are introduced to a one-line correction, reminiscent of the measurement step in Kalman Filtering, and are tasked with tuning the filter's gains to correct for this bias. While the project is deliberately short and approachable, it sets students up for more advanced topics, including ordinary differential equations (ODEs), by framing numerical integration as a stepping stone toward solving ODEs.

\subsection{Gradient Descent with Equality Constraints: From Basketball to Gymnastics and Diving} 

The second project takes students into the realm of optimization, beginning with a familiar scenario: calculating the initial speed and angle required for a basketball to make a successful free throw. Students explore multiple formulations of the problem. Initially, they solve a $2 \times 2$ system of linear equations to determine the initial speed ($v_x, v_y$) while fixing the time of flight. Next, they use gradient descent to minimize the squared distance from the basketball to the hoop, gradually adding complexity by introducing equality constraints that factor in time of flight and initial speed. 
This project then transitions to more complex optimization challenges, including modeling the dynamics of a simple gymnast as a floating-base representation of a bar with two point masses. Students optimize the time of flight and initial speed subject to linear constraints on the landing posture, utilizing the JuMP package to verify their results. The project culminates in an exciting optimization involving a hinged-model of a diver on a 10-meter platform, pushing students to refine their understanding of gradient descent with constraints. This progression from basketball to gymnastics, and finally to diving, demonstrates how optimization techniques are applied across diverse engineering problems.

\subsection{Modeling and Feedback Control of the BallBot}

The final project immerses students in the modeling and feedback control of a planar version of the BallBot, a simplified model of a robot that balances on a basketball using omni-directional wheels. Students begin by calculating the potential and kinetic energy of the robot's two components---the planar basketball and the torso---using Lagrange’s Method to derive the so-called "Robot Equations" of motion,
$$ D(q) \cdot \ddot{q} + C(q, \dot{q}) \cdot \dot{q} + G(q) = \Gamma. $$

Using software tools, students develop a linear state-variable model and transfer function representations of the BallBot. The project draws a direct comparison to the Segway, with the complete control design for the Segway provided in the course textbook. Following this methodology, students design, analyze, and test controllers for the linearized model of the BallBot, eventually extending their work to the full nonlinear model. The project offers a well-motivated introduction to feedback control, with all steps supported by Julia-based tools that make the process transparent, hands-on, and engaging.

%% file: Sections/05StudentAssessment.tex
We provide access to the raw midterm and end-of-term teaching evaluations in~\ref{sec:ROB201midtermEvaluations}
 and \ref{sec:ROB201endtermEvaluations}, respectively. The midterm assessments correspond to approximately the middle of Chapter 6, Engineering Applications of the Derivative.

\subsection{Pilot-Midterm Reviews}
Midterm teaching evaluations for the pilot offering of \textsc{Calculus for the Modern Engineer} provide initial insight into the effectiveness of the course and how well it aligns with the reform goals. With 22 of 24 students responding, the feedback highlights strong student engagement and satisfaction. The course received an overall rating of \textbf{4.29} (on a scale of 1 to 5, with 5 best), with a score of \textbf{4.50} on whether it advanced students' understanding of the subject matter. Additionally, \textbf{85\%} of students indicated that their interest in the subject had increased because of the course, reflecting its success in making calculus more relatable and relevant to modern engineering challenges.

A recurring theme in student feedback is the positive impact of the Julia-based programming assignments and projects. Several students emphasized that the integration of programming has brought calculus concepts to life by illustrating their real-world applications in robotics and engineering. One student commented, ``The multiple formats of the homework (written and Julia) were really great for understanding the math while also seeing real applications.” This fusion of theory and practice has made the course more engaging, allowing students to experience firsthand how mathematical concepts such as numerical integration and differentiation are applied in solving real engineering problems.

Students also appreciated the structured approach to the material, with one noting, ``We're never just told to 'memorize a formula'; instead, we are shown the steps needed to get there, giving a deeper understanding of various aspects of calculus.” This shift from rote memorization to a deeper, conceptual understanding is at the core of the course’s philosophy and has been well received by students. 

However, the evaluations also revealed some areas for improvement. Two students found the written homework assignments less valuable, feeling that they sometimes amounted to busywork. One student expressed a desire for more challenging coding assignments, stating that the provided starter code reduced opportunities for developing independent coding skills. These comments will inform future iterations of the course, as we seek to balance the need for structured guidance with opportunities for more advanced students to apply coding techniques more autonomously.

Overall, the midterm evaluations indicate that \textsc{Calculus for the Modern Engineer} is successfully meeting its goals of fostering engagement and advancing understanding through the integration of calculus with real-world applications. The course has been particularly effective in demonstrating how calculus can be a powerful tool in robotics and engineering when paired with modern computational tools like Julia. This feedback will be crucial in refining the course for future semesters, ensuring that it meets student expectations while preparing them for the challenges of modern engineering.


\subsection{End-of-Pilot-Term Reviews}

The final teaching evaluations for the pilot offering of \textsc{Calculus for the Modern Engineer} provide deeper insights into the course's reception and its alignment with reform goals. With 20 of 24 students responding, the feedback reflects strong engagement and overall satisfaction with the course.

The course received an overall rating of \textbf{4.8/5} for excellence, with scores of \textbf{4.7/5} for advancing students' understanding of the subject matter and \textbf{4.6/5} for providing a good understanding of concepts and principles. Importantly, \textbf{85\% of students} indicated that their interest in calculus increased as a result of the course. Students also noted that the workload was reasonable, with a median score of \textbf{2.9/5}, aligning with the median score of all courses offered in the College of Engineering.

A recurring theme in the feedback was the positive impact of the programming assignments and projects. Students emphasized that the integration of Julia programming helped connect mathematical theory to engineering applications. One student commented, \textit{``The programming assignments provided context where the math concepts would be applied and will help me with future classes and the real world.''} Another noted, \textit{``The multiple formats for the homework helped me to understand the material and know how to apply it to real work scenarios.''}

Students also appreciated the structured approach to the material. One remarked, \textit{``This course bridges theoretical concepts with practical applications. It pushed me to deepen my problem-solving and critical thinking.''} Many expressed gratitude for the textbook, praising its clarity and long-term value as a reference.

Importantly, the evaluations also revealed areas for improvement. Some students expressed concerns about the written homework assignments, finding them occasionally repetitive or less valuable compared to the programming tasks. Others suggested more challenging computational exercises to provide additional opportunities for growth. Specific comments included requests for in-class activities, quizzes, and increased participation to maintain focus during lectures.

Additional constructive criticism focused on the course's scheduling, with several students suggesting shorter and more frequent lecture sessions to improve engagement. Others recommended integrating more examples and applications earlier in the course to reinforce theoretical concepts.

Overall, the final evaluations indicate that \textsc{Calculus for the Modern Engineer} met its primary goals of fostering engagement and improving conceptual understanding through an updated, application-driven curriculum. Students particularly valued the focus on computation and real-world applications, with many recommending the course as a replacement for traditional calculus sequences.

Looking forward, this feedback highlights opportunities to:
\begin{itemize}
    \item Increase in-class participation through quizzes or group activities.
    \item Refine homework assignments to balance computational challenges with conceptual depth.
    \item Optimize the course schedule to better accommodate student attention spans and engagement.
    \item Expand programming assignments to challenge more advanced students.
\end{itemize}

The overwhelmingly positive response, combined with thoughtful critiques, will guide future iterations of the course. These insights affirm the effectiveness of this reform in calculus education, demonstrating its potential as a scalable model for integrating computational tools with mathematics instruction.

%% file: Sections/06BarriersToReform.tex
One of the most significant barriers to calculus reform is engineering's outsourcing of more than two-thirds of its first-year STEM curriculum to departments outside of engineering (e.g., mathematics, physics, chemistry, biology). Most engineering faculty accept this as normal, as it reflects their own training. Although there is logical merit in mathematics departments teaching calculus, the author’s 35 years of experience illustrate that the current arrangement leaves little incentive for faculty in either department to propose a radically different curriculum. The author had to wait until the latter stages of his research career to allocate the thousands of full-time hours needed over two summers and an academic year to develop the textbook, homework sets, and projects. Active research faculty, whose academic progress hinges on research output, cannot dedicate even a fraction of this time to teaching reform. 

Moreover, when the first-year curriculum is outsourced, its importance is often undervalued. Would significant effort devoted to curriculum improvements be recognized during performance evaluations? Without administrative support from the College of Engineering at Michigan, the author’s plans to use a sabbatical for this project would have been halted early on. Despite the pilot course's initial success, acceptance among engineering faculty outside the Robotics Department remains uncertain. Change, especially in curriculum structure, is inherently difficult.

Relying on student surveys post-graduation to assess essential mathematical skills for engineering has limited applicability to modern disciplines. These surveys often reflect the perspectives of students trained under traditional methods, making it challenging for them to envision how calculus integrates with contemporary computational tools. Furthermore, many such surveys were conducted before fields like robotics, artificial intelligence, and data science became pivotal, areas that demand an intricate blend of mathematics and programming. Consequently, these studies often miss the mark in identifying the evolving needs of the next generation of engineers.

Extensive research has analyzed the experiences of STEM students in traditional calculus courses \cite{bressoud2021calculus, hatfield2022stem}. Efforts aimed at creating a more positive learning environment and reducing the ``weed-out” culture of the traditional curriculum are essential if one believes that meaningful change can occur in teaching outcomes without fundamentally altering the course content or its instructors. Students can discern when an instructor lacks direct experience in applying calculus to practical scenarios. This distinction—between ``Talking the Talk” and ``Walking the Walk”---affects their engagement and trust in the learning process.

%% file: Sections/07Discussion.tex
One of the greatest advantages of developing \textsc{Calculus for the Modern Engineer} within the context of the University of Michigan’s newly established Robotics Department has been the freedom to break away from conventional norms in curriculum design. This freedom extended to reimagining the entire undergraduate robotics curriculum, allowing for innovative sequencing and integration of essential mathematical concepts. Notably, the emphasis on math innovation began with the launch of \textbf{ROB 101 Computational Linear Algebra} \cite{apte2021solving, grizzle_rob101, grizzle2020rob101}, two years before the department's official founding; see also~\ref{sec:ROB101}. In addition to recognizing its foundational role in robotics, artificial intelligence and the processing of sensor data, ROB 101 was created with the aims of (1) introducing linear algebra before calculus, (2) opening up mathematics to everyone because computational linear algebra builds directly on high school algebra, and (3) making math fun to learn by marrying programming with carefully curated projects.

Building on this effort, two years after the department was established, the author undertook the task of bringing similarly motivated changes to the calculus curriculum. The traditional sequence of spending three to four semesters solely on calculus is unnecessarily long, especially given the mathematical demands of modern engineering disciplines. By condensing Calculus I, II, \& IV into a one-semester intensive course, \textsc{Calculus for the Modern Engineer} allows students to allocate the freed-up semesters to exploring advanced mathematical topics. 

Importantly, the new calculus course will not reduce the overall math requirements for Robotics majors. Students are required to complete an additional 8 credits of advanced mathematics beyond \textsc{Calculus for the Modern Engineer}. While credits from Calculus I \& II are still applicable if taken before this course, pursuing these classes afterward requires advisor approval. Recommended courses for further study include Math 215 Multivariable and Vector Calculus, Math 217 Linear Algebra (with proofs), Math 312 Applied Modern Algebra, Math 351 Principles of Analysis, Math 371 (Engin 371) Numerical Methods, Math 412 Abstract Algebra, Math 416 Theory of Algorithms, and Math 451 Advanced Calculus I, which offers a proof-based perspective taught through Real Analysis in $\mathbb{R}^n$. In robotics, mathematical knowledge equates to power, especially when coupled with the ability to rapidly translate theory into code.

This autonomy in curriculum development allowed for a focus on real-world applications, computational tools, and solid theoretical foundations. Integrating topics such as numerical integration, optimization, dynamical systems, and feedback control ensures that students acquire both the theoretical and practical skills necessary for modern engineering challenges. This freedom also gave us the flexibility to explore innovative pedagogical strategies. Rather than relying on rote memorization and manual problem-solving, we adopted a project-based learning approach where students tackle real-world robotics problems from the beginning. For instance, numerical integration is introduced not as an abstract mathematical concept but as a practical tool for calculating positions and velocities from acceleration data in robotic systems. This focus on practical, goal-oriented learning makes calculus more relevant and accessible to engineering students.

Central to the new calculus course is the integration of computation, specifically through the Julia programming language. Julia’s high-performance capabilities allow students to visualize and manipulate complex mathematical expressions, akin to how they would work on paper, which would be difficult to comprehend through traditional pen-and-paper methods. As an open-source and functional programming language, Julia ensures that expressions in code closely resemble their mathematical counterparts, making it easier for students to translate mathematical ideas into computational implementations. Students apply calculus concepts directly to tasks like optimizing multivariable functions, deriving differential equations for mobile robots, and developing feedback controllers—skills essential for robotics and modern engineering fields.

Unlike other recent efforts on integrating computation into the traditional calculus curriculum, such as \textit{Calculus with Julia} and the {\it MTH 229: Calculus Computer Laboratory}, \textsc{Calculus for the Modern Engineer} goes further by fundamentally reordering the curriculum to meet the needs of modern engineering students. For example, the course begins with definite integration, a topic that resonates with students due to its intuitive connection to sums and areas under curves. This integration-first approach provides a solid foundation before introducing more abstract concepts such as limits, differentiation, and differential equations.

To ensure students are well-rounded, the course strikes a balance between computational fluency and manual problem-solving. While students gain extensive experience with computational tools like Julia, they also engage in enough hand computation to succeed in more traditional engineering courses. This dual focus prepares students not only for academic success, but also for the engineering challenges they will face in industry and research.

Ultimately, \textsc{Calculus for the Modern Engineer} aims to transform calculus from a conceptual hurdle into an empowering tool that is exciting and interesting to learn. Through a combination of Julia-based programming, robotics case studies, and coding-based projects, students leave the course with a deep conceptual understanding of calculus and the ability to apply it meaningfully in their future careers.

%% file: Sections/08ConclusionPerspectives.tex
\textsc{Calculus for the Modern Engineer} seeks to align calculus education with the computational and practical needs of contemporary engineering. By rethinking the traditional structure and emphasizing project-based learning, computation, and real-world applications, this course offers a blueprint for modernizing the mathematics curriculum for engineers. The successful pilot has shown that integrating computational tools like Julia with rigorous mathematical principles not only deepens student understanding but also enhances their ability to apply calculus meaningfully in engineering contexts.

To build on this success, a GitHub repository with a common open-source license will be launched, providing access to the textbook source files, written homework sets, and Jupyter notebooks for programming assignments and projects (the link will be posted here~\cite{grizzle2025rob201} when available) This initiative aims to broaden the course's reach and invite contributions from educators and the global learning community. 

Future iterations of the course will be refined through student feedback, such as additional coding challenges for advanced learners and expanding case studies and projects to include applications in aerospace, mechanical engineering, or biomedical fields to further showcase the broad applicability of calculus across various disciplines.

In the second offering of the pilot course (Winter 2025), we will introduce a low-stakes, timed final exam worth 10\% of the final grade. The purpose of this exam is to assess how well students can recall and apply the knowledge gained through homework assignments and projects. The unfortunate need for this addition has emerged through multiple channels. Some students in the earlier ROB 101 course reported that it felt ``content free'' when discussing with advisors in CSE. While the actual content is substantial (see \ref{sec:ROB101}), undergraduate teaching assistants in ROB 201 have suggested that these comments likely mean the students felt they could complete the course without fully engaging with the material. This needs to be investigated. The final exam’s format is not yet finalized, but it will likely be electronic, allowing students access to a Google Doc containing their own notes.

Reforming calculus is challenging, not least because institutional support can be difficult to secure. A senior math colleague, committed to reform, candidly expressed initial doubts, predicting that the course’s compressed format and the diversity of students’ preparation levels might lead to a 'trainwreck.' Concerns ranged from insufficient depth and lack of foundational calculus to over-reliance on programming tools. Yet, the final evaluations decisively disproved these fears. Students overwhelmingly praised the course’s structure, rigor, and modern approach. As one student put it, \textit{``The course was challenging yet rewarding, pushing me to deepen my problem-solving skills and critical thinking when it comes to calculus.''} 

In summary, \textsc{Calculus for the Modern Engineer} represents a potentially transformative approach to teaching calculus, equipping students with computational proficiency and an appreciation for how mathematical theory underpins practical problem-solving. Michigan Robotics is committed to continually evolving the course, fostering an educational experience that prepares students to excel in the complex landscape of modern engineering.

%% file: Sections/09ContentRationnale.tex

The process of selecting content for \textsc{Calculus for the Modern Engineer} was not a simple exercise in taking the traditional list of calculus topics, crossing out the unnecessary ones, and adding a few applications. Such an approach would be akin to subtractive manufacturing—like milling away excess material—leaving behind a functional yet disjointed curriculum. Instead, an additive approach was taken, similar to 3D printing, where the curriculum was built layer by layer, with clear milestones and goals guiding the structure. A key feature of this design was the decision not to adhere to the traditional order of calculus topics but rather to prioritize content that would best serve modern engineering students, especially in robotics.

The textbook begins with a rigorous treatment of pre-calculus concepts, presented in greater depth than what is typically encountered at the high school level. This strong foundation ensures that all students, regardless of background, start with a clear grasp of fundamental mathematical principles such as notation, functions, and algebra. Historical context is interwoven into the material to enhance students' understanding. For instance, Bernoulli's discovery of $e$ and Archimedes' method of approximating $\pi$ are introduced not merely as mathematical curiosities but as milestones in the historical development of ideas central to modern calculus. Archimedes' technique---using inscribed and circumscribed polygons to provide bounds for $\pi$---serves as a precursor to modern concepts like Riemann sums and the Squeeze Theorem. By framing these historical insights as bridges to modern ideas, the course makes abstract calculus concepts more accessible and intuitive.

One of the foundational innovations of the course is the decision to begin with integration instead of differentiation, because sums are easy to understand. Moreover, they rely on limits at infinity, which are easier to understand than singe-sided limits at points, especially in the context of rational functions, which are obtained when performing Riemann sums for monomials and the exponential function. Limits at infinity are taught using both the formal epsilon-delta method and a more intuitive means by inspection. These two approaches, humorously dubbed the ``hard way" and the ``easy way," demystify rigorous proofs and help students understand when each method is most appropriate. The second chapters's focus on limits at infinity, which describe the long-term behavior of functions, therefore lays the groundwork for Chapter 3 on Riemann sums and prepares students for more complex topics in calculus.

After developing integration via Riemann-Darboux sums, more efficient techniques for numerical integration are explored in Julia. This is followed by concrete applications of the Riemann integral to subjects such as center of mass, moments of inertia, and volumes of revolution. Through these applications, students became acquainted with integration as a rigorous (and painless) way to add up differential quantities to form a whole. Formulas such as 
$$
    I_z := \int_a^b dI_z =  \rho \cdot h\cdot \int_a^b  \left[x^2 \cdot \left( f(x) - g(x)\right)   + \frac{1}{3} \cdot \left( f^3(x) - g^3(x)\right) \right] \,  dx
$$
for the moment of inertia about the $z$-axis of a uniform planar body upper bounded  $f:[a, b] \to \real$ and lower-bounded $g:[a, b] \to \real$ cannot be reasonably computed by hand for interesting shapes. However, in code, they become manageable and meaningful. Because the computation of $dI_z$ for a rectangle of height $f(x)-g(x)$ and infinitesimal width $dx$ involves a first integration of ``point masses'' along the $y$-axis, students are offered a peek at what a course on multivariable calculus addresses. 

It is important to note that at this point, students know how to numerically compute a definite integral of just about any function, but they only know how to compute the \textbf{indefinite integral} of polynomials and exponentials. In an optional-read section, the indefinite integral of sine and cosine are illustrated via Euler's formula. The concepts of antiderivatives and the fundamental theorems of calculus are delayed until Chapter 7, in part out of necessity because derivatives have not yet been introduced, and in part because students coming out of a standard calculus curriculum wrongly conflate integration with antiderivatives, leading them to shy away from hard integrals such as the one shown above for the moment of inertia. \textsc{Calculus for the Modern Engineer} seeks to avoid this pedagogical pitfall by making a clear delineation between the two notions of integration. 

\emstat{Intuitive, heuristic explanations are provided for each of the topics in the textbook. Students who do not wish to understand why calculus ``works as it does'' can still understand  ``how calculus works'' through the projects, the Julia programming assignments, and the written homework. So far in the pilot term, a majority of students have been grateful for the insight provided by the more fundamental approach. Some have even commented that the textbook will grow with them as they mature and progress through their engineering studies.}

After ``numerical'' integration, the mid-course goal is to introduce students to optimization techniques because it is a core tool in modern engineering and robotics. Hence, the focus of Chapter 4 is on the tools required for differentiation of functions of a single variable. Single-sided limits are  first introduced using the intuition students have acquired with limits at infinity, namely
the limit from the right is the same as the following (positive) infinite limit
    $$ \lim_{x \to x_0^+} f(x) := \lim_{\eta \to + \infty} f(x_0 + \frac{1}{\eta}),$$
    and the limit from the left is the same as the corresponding negative infinite limit
    $$ \lim_{x \to x_0^-} f(x) := \lim_{\eta \to -\infty} f(x_0 + \frac{1}{\eta}).$$
    The key points being
    \begin{itemize}
        \item $\frac{1}{\eta}$ is never exactly equal to zero, so one never evaluates the function at $x_0$;
        \item  as $\eta \to \infty$,  $x_0 + \frac{1}{\eta}$ approaches $x_0$ from the \textbf{right} because $x_0 + \frac{1}{\eta}> x_0$ for all $\eta >0$; and
        \item  as $\eta \to -\infty$,  $x_0 + \frac{1}{\eta}$ approaches $x_0$ from the \textbf{left} because $x_0 + \frac{1}{\eta}><x_0$ for all $\eta <0$.
    \end{itemize}
Later the substitution ${ h:=\frac{1}{\eta} }$ is made, so that ${ \eta \to +\infty \iff h \to 0^+ }$ (${h}$ decreases to zero through positive values, that is, ${h}$ approaches zero from the right) and ${ \eta \to -\infty \iff h \to 0^- }$ (because ${h}$ increases to zero through negative values, that is, ${h}$ approaches zero from the left), leading to the classical forms of the definitions. Continuity of functions at a point is treated carefully enough that limits can be taken inside of continuous functions. With this tool at hand, limits that will be needed in a chapter on differentiation are addressed. Without a doubt, Chapter 4 is the most challenging chapter in the textbook. This is clearly acknowledged with warnings in writing and verbally in lecture! A key point however, is that this highly technical material is delayed until students have become accustomed to using computation to augment their ability to break concepts into manageable pieces. Specifically, the standard epsilon-delta notions of limits are heuristically illustrated in code by adding a finite array of ``vanishing'' positive and negative numbers to a given point, and evaluating a function on the two arrays. 

With this foundation, Chapter 5 treats differentiation of single-variable functions from the perspective of rise over run and linear approximations. After working the obligatory examples using the definition, the Rules of Differentiation are introduced. Derivations for the product rule and the chain rule are illustrated using the linear approximation of a function about a point. Between the definition and the rules, the students have Table~\ref{tab:CommonFunctionsAndTheirDerivatives}, which represents about the only thing they are asked to memorize in the course.
\renewcommand{\arraystretch}{1.5}
\begin{table}[htb]
\centering
\begin{tabular}{|c |c|c|c|c|c|c|c|c|}
\hline
Function & $c$ & $x^n$ & $e^x$ & $\ln(x)$ & $\sin(x)$ & $\cos(x)$ & $\tan(x)$ & ${\rm atan}(x)$\\
\hline
Derivative & $0$ & $nx^{n - 1}$ & $e^x$ & $\frac{1}{x}$ & $\cos(x)$ & $-\sin(x)$ & $1 + \tan^2(x)$ & $\frac{1}{1 + x^2}$\\
\hline
\end{tabular}
\caption{Common Elementary Functions and their Derivatives (worth committing to memory to showcase your expertise).}
\label{tab:CommonFunctionsAndTheirDerivatives}
\end{table}
Using the Symbolics package in Julia, students could differentiate essentially ``anything''.

Applications of single-variable derivatives include extreme points of functions and arc or path length. The formula for path length is derived from first principles and specialized to the standard cases. Its integral, 
$$ S := \int_{t_0}^{t_f}  \sqrt{\left(\frac{dx}{dt}\right)^2 + \left(\frac{dy}{dt}\right)^2} \, dt $$
is very hard to evaluate by hand for anything more complicated than a circular arc, but can be done numerically with three lines of code in Julia. Its use to calculate the length of a path followed by a mobile robot is easy to motivate. In addition, students are pleasantly surprised to see the Pythagorean Theorem show up in such a clever manner.

The primary goal in the study of optimization is gradient descent with equality constraints. This requires that partial derivatives be treated after derivatives of single-variable functions. With ROB 101 \textit{Computational Linear Algebra} as a prerequisite, the connection to single variable derivatives is straightforward, namely, 
$$
\begin{aligned}
\frac{\partial f(x_0)}{\partial x_i} :=& \lim_{h \to 0} \frac{f(x_0 + h\cdot e_i) - f(x_0)}{h} \\
=& \frac{d}{d h} f(x_0 + h\cdot {e_i}) \big|_{h = 0)}
\end{aligned}
$$
assuming the limit exists and is finite; here, $e_i$ is the natural basis for $\mathbb{R}^n$. 

\emstat{In the pilot semester, a first-semester student from the Ross School of Business took the course without ROB 101 or high school linear algebra, but with AP Calclulus A and a high school programming class. They initially struggled to master Julia, but by the time differentiation was treated, they were fine. The teaching assistants in the course spent extra time with the student to explain vectors and matrices. At the time of the writing of this article, they are near the top of the class. }

Chapter 6 is devoted to engineering applications of the derivative. Root finding and gradient descent are studied and illustrated on concrete multivariable examples. Then, constrained optimization via Lagrange multipliers and gradient descent with equality constraints are studied and applied to meaningful examples. As an application of the total derivative in combination with partial derivatives, a very simplified introduction is made to Lagrange's method for deriving equations of motion. Its implementation in the Symbolic package of Julia is emphasized, along with special purpose tools to create write out the differential equations as a function in Julia. At this points, ODEs are equations that contain derivatives, with illustrations coming from $F=ma$ in the form $F= m \dot{v}$ or $F = m\ddot{x}$, in addition to the awe-inspiring robot equations. To be very clear, this topic cannot be treated without the appropriate software tools because the hand calculations are just too arduous for any interesting examples.

Chapter 7 treats antiderivatives and the fundamental theorems of calculus. Students are shown how to use LLMs, Wolfram Alpha Pro, and Julia to compute antiderivatives, and then they are shown basic cases of u-substitution, integration by parts, trig-substitution, and partial fraction expansion for real roots. The goal here is to provide adequate fluency for understanding derivations in their follow-on engineering courses, with a wink-wink to the fact that one day, these methods might suffer the same fate as the computation of square roots by hand. But until that day arrives, they need to compete with students coming through a standard three-semester calculus sequence. Chapter 8 treats improper integrals, with the motivating application coming from probability.

The course culminates in Chapters 9 and 10, treating ODEs, and Laplace transforms through the lens of feedback control, feeding into the third and final project, modeling and control of the BallBot, a robot they study in \textit{ROB 311: Build Robots and Make Them Move}. This actually solves a challenge in ROB 311, where students do not have adequate time to build a principled feedback controller for their robot....it remained wobbly due to poorly understood concepts in PD control.

This described structure provided an overall framework to work within. As the course was fleshed out, it was not uncommon to uncover gaps or missing material---topics that were necessary to set the stage for upcoming concepts. This led to an iterative process where earlier chapters were revised and expanded to better prepare students for the later material. Each addition was intentional, adding depth and context to ensure that the flow of learning remained coherent and that students would fully understand how the concepts build upon one another. This iterative approach allowed the course's content to evolve organically, ensuring that the final product was both lean, innovative, and aligned with the real-world needs of robotics students. 

With the topics building on one another, fewer tangents were introduced simply because they were part of the traditional curriculum. However, as mentioned above, it was also important to ensure that students could still compete effectively with their peers who had followed a more conventional calculus path. This meant incorporating traditional manual calculation methods for differentiation and integration of common, elementary functions. The written homework assignments reflect this need, giving students practice with these foundational techniques, while the programming assignments focus on more realistic, complex examples and employ software tools to handle the heavier computations.

In addition to the careful selection of topics, the course emphasizes a balance between theoretical rigor, intuition, and practical application. Students not only learn the mathematical underpinnings of the formulas and techniques they use, but also how to implement them in code, making the abstract more concrete. This dual focus allows students to gain a deeper understanding of the material while developing practical skills that they can carry into their professional lives. Whether they are solving a system of ODEs or analyzing sensor data from a drone, students leave this course with the confidence that they can apply their calculus knowledge to a wide range of real-world problems.

Ultimately, the content of \textsc{Calculus for the Modern Engineer} reflects an intentional restructuring of calculus education. By focusing on the specific needs of robotics students and emphasizing computational tools, the course equips students with the skills and knowledge they need to thrive in a technology-driven world. This approach ensures that students see calculus not as a series of disconnected, theoretical concepts but as an integrated, powerful tool for solving the problems they will encounter as engineers.

The author believes that the selection of topics is appropriate for students in Aerospace Engineering, Mechanical Engineering, and Naval Architecture and Marine Engineering as well as Robotics. Because the source material is being placed in an open-source repository, an enterprising colleague can replace mechanically-oriented material with examples of their choice. If the collaboration with the Ross Business School grows, the author will work with their faculty to expand the set of examples.

%% file: Sections/10ROB101ComputationalLinearAlegbra.tex
Computational Linear Algebra is a first-semester, first-year undergraduate course that demonstrates the combined power of mathematics and computation for reasoning about data and making discoveries about the world. Linear algebra and coding are rapidly becoming essential foundations for the modern engineer in a computational world \cite{grizzle_rob101, grizzle2020rob101}. Students in this course gain insights into the mathematical theory of linear algebra and its realization in practical computational tools. While mathematics is the language of engineering, Prof. Chad Jenkins of Michigan Robotics underlines that coding is believing and realizing it. 
The mathematical content of ROB 101 revolves around systems of linear equations, their representation as matrices, and numerical methods for their analysis. These methods come to life through the lens of robotics and contemporary intelligent systems, showcasing compelling applications.

The course design embodies inclusivity, with a curriculum that assumes a background in Algebra but none in programming. It opens the doors of Linear Algebra and computational tools to almost anyone in a STEM field at U-M. By employing Julia, a browser-friendly programming language, all students, regardless of access to high-end computing resources, could participate effectively. 
In partnership, a third of the Fall 2020 pilot class participated remotely from Morehouse and Spelman Colleges. The course features three major projects, connecting the theoretical foundation of linear algebra to applications in fields such as computer vision, climate data analysis, and feedback control systems.

This approach not only prepares students for advanced topics, but also reimagines how mathematics is introduced to first-semester engineering undergraduates. ROB 101 embodies a commitment to equity, innovation, and excellence in teaching.

\subsection{Contents Overview}
The ROB 101 curriculum is structured to introduce first-year students to computational linear algebra as a foundational tool for robotics and engineering. Its chapters are organized as follows:

\begin{enumerate}

\item {\bf Introduction to Systems of Linear Equations}  
\textbf{Notes:} $\bullet$ 2 hours lecture + recitation. $\bullet$ Introduces foundational algebraic concepts for systems of equations. $\bullet$ Examples cover unique, no, and infinite solutions. $\bullet$ Includes a programming exercise for basic matrix operations in Julia.  
\textbf{Learning Objectives:} $\bullet$ Recognize and classify systems of linear equations. $\bullet$ Develop solution methods, including substitution and elimination. $\bullet$ Appreciate scalability limits of hand-methods. $\bullet$ Build familiarity with applications in robotics and engineering.  
\textbf{Content:} Unique, no, and infinite solutions; substitution, elimination, graphical interpretations, and $Ax = b$ as a framework. (Optional Reads: Historical relevance in robotics.)  

\item {\bf Vectors, Matrices, and Determinants}  
\textbf{Notes:} $\bullet$ 4 hours lecture + recitation. $\bullet$ Matrix representation and operations for systems of equations. $\bullet$ Focus on determinants as a key property.  
\textbf{Learning Objectives:} $\bullet$ Understand vectors and matrices as data structures for linear systems. $\bullet$ Explore determinant properties. $\bullet$ Represent and solve systems in matrix form.  
\textbf{Content:} Vectors, matrices, operations, determinants, and solvability. (Optional Reads: Additional determinant properties.)  

\item {\bf Triangular Systems of Equations: Forward and Back Substitution}  
\textbf{Notes:} $\bullet$ 4 hours lecture + recitation. $\bullet$ Solving triangular systems as simplified linear systems. $\bullet$ Forward and back substitution implemented in Julia.  
\textbf{Learning Objectives:} $\bullet$ Solve lower and upper triangular systems. $\bullet$ Develop substitution methods. $\bullet$ Compute determinants by inspection. $\bullet$ Understand triangular matrices in numerical methods.  
\textbf{Content:} Forward and back substitution, triangular matrices, and applications. (Optional Reads: Numerical stability.)  

\item {\bf Matrix Multiplication}  
\textbf{Notes:} $\bullet$ 3 hours lecture + recitation. $\bullet$ Mechanics and significance of matrix multiplication. $\bullet$ Applications include transformations and compositions.  
\textbf{Learning Objectives:} $\bullet$ Perform standard matrix multiplication. $\bullet$ Understand row-column relationships. $\bullet$ Explore outer-product formulation for LU factorization. $\bullet$ Apply matrix multiplication to practical problems.  
\textbf{Content:} Matrix multiplication views, examples, and applications. (Optional Reads: Outer-product formulation proof.)  

\item {\bf LU (Lower-Upper) Factorization}  
\textbf{Notes:} $\bullet$ 3 hours lecture + recitation. $\bullet$ LU factorization without pivoting; pivoting explained using Julia tools.  
\textbf{Learning Objectives:} $\bullet$ Understand LU factorization. $\bullet$ Solve systems using LU decomposition. $\bullet$ Relate LU factorization to triangular systems.  
\textbf{Content:} LU factorization, solving $Ax = b$, and numerical applications. (Optional Reads: Pivoting algorithms.)  

\item {\bf Determinants, Matrix Inverses, and Transposes}  
\textbf{Notes:} $\bullet$ 2 hours lecture + recitation. $\bullet$ Matrix properties and implications for linear systems.  
\textbf{Learning Objectives:} $\bullet$ Compute determinants and inverses. $\bullet$ Apply transposes. $\bullet$ Recognize these properties’ roles in solutions.  
\textbf{Content:} Determinants, inverses, transposes, and computations. (Optional Reads: Matrix identities.)

\item {\bf The Vector Space $\mathbb{R}^n$: Part 1}  
\textbf{Notes:} $\bullet$ 4 hours lecture + recitation. $\bullet$ Introduces students to the abstract properties of vector spaces.  
\textbf{Learning Objectives:} $\bullet$ Define and understand vector spaces. $\bullet$ Analyze linear independence and span. $\bullet$ Relate vector spaces to solutions of $Ax = b$.  
\textbf{Content:} Definition of vector spaces, linear independence, span, and applications in robotics. (Optional Reads: Basis and dimension.)  

\item {\bf Norms and Least-Squares Solutions}  
\textbf{Notes:} $\bullet$ 3 hours lecture + recitation. $\bullet$ Introduces optimization concepts through least-squares methods.  
\textbf{Learning Objectives:} $\bullet$ Compute norms and interpret their significance. $\bullet$ Apply least-squares methods to solve overdetermined systems. $\bullet$ Explore regression as a practical application.  
\textbf{Content:} Vector norms, least-squares solutions, regression, and applications. (Optional Reads: Generalized least-squares methods.)  

\item {\bf The Vector Space $\mathbb{R}^n$: Part 2}  
\textbf{Notes:} $\bullet$ 5 hours lecture + recitation. $\bullet$ Explores orthogonality and its applications in linear algebra.  
\textbf{Learning Objectives:} $\bullet$ Understand orthogonal vectors and matrices. $\bullet$ Apply the Gram-Schmidt process. $\bullet$ Explore QR factorization.  
\textbf{Content:} Orthogonality, Gram-Schmidt process, QR factorization, and practical applications. (Optional Reads: Modified Gram-Schmidt.)  

\item {\bf Eigenvalues, Eigenvectors, and Advanced Concepts}  
\textbf{Notes:} $\bullet$ 3 hours lecture + recitation. $\bullet$ Introduces eigenvalues and eigenvectors.  
\textbf{Learning Objectives:} $\bullet$ Compute eigenvalues and eigenvectors. $\bullet$ Understand their role in linear transformations. $\bullet$ Apply these concepts in robotics and dynamics.  
\textbf{Content:} Eigenvalues, eigenvectors, diagonalization, and applications. (Optional Reads: Spectral theorem.)  

\item {\bf Solutions of Nonlinear Equations}  
\textbf{Notes:} $\bullet$ 3 hours lecture + recitation. $\bullet$ Introduces methods for solving nonlinear systems. $\bullet$ Shows how linear approximations of nonlinear functions extend linear algebra techniques to broader problems.  
\textbf{Learning Objectives:} $\bullet$ Introduce numerical differentiation and linear approximations of nonlinear functions. $\bullet$ Explore numerical methods for nonlinear equations. $\bullet$ Apply Newton-Raphson method. $\bullet$ Understand Jacobians and gradients through numerical differentiation (central or symmetric differences).  
\textbf{Content:} Bisection, Newton-Raphson, Jacobians, and gradients. (Optional Reads: Applications in robotics, such as calibration of a LiDAR sensor.)  

\item {\bf Basics of Optimization}  
\textbf{Notes:} $\bullet$ 4 hours lecture + recitation. $\bullet$ Introduces optimization methods in robotics.  
\textbf{Learning Objectives:} $\bullet$ Apply gradient descent. $\bullet$ Understand the role of the Hessian matrix. $\bullet$ Solve constrained and unconstrained optimization problems using tools in Julia.  
\textbf{Content:} Gradient descent, Hessians, optimization in robotics. (Optional Reads: Quadratic programming.)  

\item {\bf Background for Machine Learning}  
\textbf{Notes:} $\bullet$ 2 hours lecture + recitation. $\bullet$ Provides foundational concepts for machine learning.  
\textbf{Learning Objectives:} $\bullet$ Understand separating hyperplanes and classifiers. $\bullet$ Understand signed-distance. $\bullet$ Apply max-margin classifiers. $\bullet$ Explore orthogonal projections.  
\textbf{Content:} Separating hyperplanes, max-margin classifiers, orthogonal projection, and applications. (Optional Reads: Soft-margin classifiers.)  

\end{enumerate}

\subsection{Three Pillar Projects in ROB 101}

ROB 101 bridges the gap between theoretical linear algebra and its applications through three cornerstone projects, each designed to highlight the real-world relevance of matrix mathematics. These projects immerse students in practical problem-solving, computational implementation, and a taste of robotics and intelligent systems.

\textbf{Building 3D Maps from LiDAR Data:} The first project, focusing on computer vision and mapping, introduces students to the processing of LiDAR data. Through this project, students learn to generate and manipulate point clouds, a foundational concept in robotics and autonomous systems. Using tools in Julia, they transform raw LiDAR scans into structured maps using matrix-vector multiplication, gaining hands-on experience with matrix transformations and 3D visualizations. This project not only reinforces linear algebra concepts but also shows how these tools power technologies such as autonomous vehicles and drones.

\textbf{Regression of Large Data Sets:} The second project shifts attention to machine learning and climate data analysis. Here, students work with real-world climate datasets to apply regression techniques based on least-squares optimization. By exploring the mathematical underpinnings of these methods, students discover the power of linear algebra in making predictions and uncovering patterns in large datasets. This project serves as a gateway to understanding the role of data-driven decision-making in engineering and science, laying the groundwork for advanced study in machine learning and artificial intelligence.

\textbf{Stabilizing a Simple Model of a Segway via MPC:} The third and final project is a capstone experience that challenges students to design a feedback control system for a Segway. Students begin by modeling the dynamics of a Segway, applying concepts of state-variable representation and discrete dynamics. They then use model predictive control (MPC) techniques, framed as underdetermined least-squares optimization problems, to develop a controller capable of balancing and stabilizing the system. This project culminates in the implementation of the control algorithm in Julia, providing students with a tangible understanding of how mathematics, coding, and engineering come together to solve complex problems.

Through these projects, students not only deepen their understanding of linear algebra but also experience the excitement of applying mathematical theory to modern engineering challenges. Each project is carefully designed to inspire curiosity, build confidence, and provide a glimpse into the possibilities awaiting them in robotics and intelligent systems.

\subsection{Course Completion Rate}

The completion rate for ROB 101 is approximately 98\% after the week-3 drop-date when students are still exploring courses, reflecting the course's design to support student success. Students who do not finish the course typically fall into two main categories: those who, despite extensive counseling, overcommit by enrolling in too many credits, and those who face unforeseen family or health challenges.

To ensure robust support, the course is staffed by a dedicated team of undergraduate student assistants who provide extensive office hours. These include late-night sessions from 8 PM to midnight on the three evenings before any assignment is due, alongside standard office hours throughout the week. This structured and accessible support system has been instrumental in maintaining the high completion rate and fostering a collaborative and encouraging learning environment.

\subsection{How the Course Came About}

The idea for this course originated with Prof. Chad Jenkins, with encouragement from then Associate Dean for Undergraduate Education, Prof. Joanna Millunchick. The textbook was composed by Prof. Jessy Grizzle. The HWs and projects were developed collaboratively with Asst. Prof. Maani Ghaffari, rising sophomore Kira Biener, MS students Madhav Achar, Fangtong Liu, and Shaoxiong Yao, and PhD students Tribhi Kathuria, Eva Mungai, Grant Gibson, and Wami Ogunbi. The course was piloted in Fall 2020 in collaboration with Prof. Dwayne Joseph, Morehouse College. Fourteen amazing students from Morehouse and Spelman joined via Zoom 38 first-year students from Michigan to comprise the pilot class. For approximately a dozen of the Michigan students in Fall 2020, this was the only course they had with an in-person section. Under the expert guidance of Lecturer III Jamie Budd and Assoc. Prof. Ram Vasudevan, ROB 101 now welcomes over 450 Michigan students each academic year and around 70 students from Howard University, Morehouse College, and Florida Agricultural and Mechanical University (FAMU). 